\documentclass{amsart}

\usepackage[german, english]{babel}
\usepackage{amssymb,amsmath,amsthm,amscd}

\input{diagrams}

\advance\oddsidemargin by -0.5cm
\advance\evensidemargin by -0.5cm
\newtheorem{theorem}{Theorem}
\newtheorem{proposition}{Proposition}
\newtheorem{corollary}{Corollary}
\newtheorem{lemma}{Lemma}
\newtheorem{definition}{Definition}
\newtheorem{assumption}{Assumption}
\theoremstyle{definition}
\newtheorem{remark}{Remark}
\newtheorem{example}{Example}

\newcommand{\modulus}[1]{\left\lvert #1 \right\rvert}
\newcommand{\norm}[1]{\left\| #1 \right\|}
\newcommand{\mklm}[1]{\left\{ #1 \right\}}
\newcommand{\eklm}[1]{\left\langle #1 \right\rangle}

\renewcommand{\d}{\,d}
\newcommand{\N}{{\mathbb N}}
\newcommand{\Z}{{\mathbb Z}}
\newcommand{\C}{{\mathbb C}}
\newcommand{\R}{{\mathbb R}}

\newcommand{\B}{{\mathcal B}}
\newcommand{\D}{{\mathcal D}}
\newcommand{\E}{{\mathcal E}}

\renewcommand{\O}{{\mathcal O}}

\newcommand{\1}{{\bf 1}}
\renewcommand{\epsilon}{\varepsilon}
\renewcommand{\phi}{\varphi}
\renewcommand{\rho}{\varrho}

\newcommand{\Cinft}{{\rm C^{\infty}}}
\newcommand{\CT}{{\rm C^{\infty}_c}}
\newcommand{\Cinftv}{{\rm \dot C^{\infty}}}
\newcommand{\CTv}{{\rm \dot C^{\infty}_c}}
\newcommand{\Cvan}{{\rm C_0}}

\renewcommand{\L}{{\rm L}}

\newcommand{\Lloc}{{\rm L^1_{\text{loc}}}}
\renewcommand{\S}{{\mathcal S}}
\newcommand{\K}{{\mathcal K}}
\newcommand{\Sym}{{\rm S}}
\newcommand{\Syms}{{\rm S^{-\infty}}}
\newcommand{\Symsl}{{\rm S^{-\infty}_{la}}}

\newcommand{\GL}{\mathrm{GL}}
\newcommand{\SL}{\mathrm{SL}}
\newcommand{\SO}{\mathrm{SO}}
\renewcommand{\O}{{\mathrm O}}

\newcommand{\g}{{\bf \mathfrak g}}
\renewcommand{\k}{{\bf \mathfrak k}}
\newcommand{\p}{{\bf \mathfrak p}}
\newcommand{\U}{{\mathfrak U}}

\newcommand{\Ad}{\mathrm{Ad}\,}

\renewcommand{\det}{\mathrm{det}\,}

\renewcommand{\Re}{\mathrm{Re}\,}

\DeclareMathOperator{\supp}{supp}

\DeclareMathOperator{\gd}{\partial}

\DeclareMathOperator{\grad}{grad}

\newcommand{\e}[1]{\,{\mathrm e}^{#1}\,}
\newcommand{\dbar}{{\,\raisebox{-.1ex}{\={}}\!\!\!\!d}}

\begin{document}

\author{Pablo Ramacher}
\title{Pseudodifferential operators on  prehomogeneous vector spaces}
\address{Pablo Ramacher, Humboldt--Universit\"at zu Berlin, Institut f\"ur Reine Mathematik, Sitz: Ru\-dower Chaussee 25, D--10099 Berlin, Germany}
\email{ramacher@mathematik.hu-berlin.de}
\thanks{The author wishes to thank Thomas Friedrich and Werner Hoffmann for many valuable discussions. Supported by the SFB 288 of the DFG}

\begin{abstract}  Let $G_\C$ be a connected, linear algebraic group defined over $\R$, acting regularly on a  finite dimensional vector space $V_\C$ over $\C$ with $\R$-structure $V_\R$. Assume that $V_\C$ posseses a  Zariski-dense orbit, so that $(G_\C,\rho,V_\C)$ becomes a prehomogeneous vector space over $\R$. We consider the left regular representation $\pi$ of the group of $\R$-rational points $G_\R$ on the Banach space $\Cvan(V_\R)$ of continuous functions on $V_\R$ vanishing  at infinity, and study the convolution operators $\pi(f)$, where $f$ is a rapidly decreasing function on the identity component of $G_\R$. Denote the complement of the dense orbit by $S_\C$, and put $S_\R=S_\C\cap V_\R$. It turns out that the restriction of $\pi(f)$ to $V_\R-S_\R$ is a smooth operator. Furthermore, if $G_\C$ is reductive, and $S_\C$ and $S_\R$ are irreducible hypersurfaces, $\pi(f)$ corresponds, on each connected component of $V_\R-S_\R$, to a   totally characteristic pseudodifferential operator. We then investigate the restriction of the  Schwartz kernel of $\pi(f)$ to the diagonal. It defines   a distribution on $V_\R-S_\R$ given by some power $|p(m)|^s$ of a relative invariant $p(m)$ of  $(G_\C,\rho,V_\C)$ and, as a consequence of the fundamental theorem of prehomogeneous vector spaces, its extension to $V_\R$,  and the complex $s$-plane, satisfies functional equations. A trace of $\pi(f)$ can then be defined by subtracting the singular contributions of the poles of the meromorphic extension.

\end{abstract}

\maketitle

\tableofcontents

\section{Introduction}

Let $G_\C$ be a connected, linear algebraic group defined over $\R$, $V_\C$ a $n$-dimensional vector space over $\C$ with $\R$-structure $V_\R$, and $\rho: G_\C \rightarrow \GL(V_\C)$ a $\R$-rational representation of $G_\C$ on $V_\C$ with a Zariski-dense $G_\C$-orbit.  The triple $(G_\C, \rho,V_\C)$ then represents a prehomogeneous vector space, and we denote the dual prehomogeneous vector space by $(G_\C,\rho^\ast, V_\C^\ast)$. In case that $G_\C$ is reductive, and the singular set $S_\C$ is an irreducible hypersurface, there exists an irreducible, homogeneous polynomial $p$ such that $S_\C=\mklm{m \in V_\C:p(m)=0}$, and $p(V_\R) \subset \R$. The singular set $S_\C^\ast$ of the dual prehomogeneous vector space is then also an irreducible hypersurface given as the set of zeros of a homogeneous, irreducible polynomial $p^\ast$ satisfying $p^\ast(V_\R^\ast)\subset \R$. The rational functions $p$ and $p^\ast$ are called relative invariants. 
Set $S_\R=V_\R\cap S_\C$,  $S^\ast_\R=V^\ast_\R\cap S^\ast_\C$, denote by $G_\R$ the group of $\R$-rational points of $G_\C$, and let $G\subset G_\R$ be the connected component containing the unit element. $V_\R-S_\R$ and  $V^\ast_\R-S^\ast_\R$ decompose into the same  number of connected components denoted by $V_i$, respectively $V_i^\ast$, $i=1,\dots,l$, each of them being a $G$-orbit. The fundamental theorem of prehomogeneous vector spaces \cite{kimura}, regarding the Fourier transform of a complex power of a relative invariant, then states that for rapidly decreasing functions $\phi$, $\phi^\ast$ on $V_\R$, respectively $V_\R^\ast$, the integrals
\begin{equation}
\label{intr1}
  \int_{V_i} |p(m)|^s \phi(m) \d m, \qquad \int_{V^\ast_i} |p^\ast(\xi)|^s \phi^\ast(\xi) \d \xi,
\end{equation}
which converge for $\Re s >0$, can be extended analytically to meromorphic functions on the whole complex $s$-plane, and satisfy the functional equations
\begin{equation}
\label{intr2}
  \int _{V_i^\ast} |p^\ast(\xi)|^{s-\frac n {\deg p}} \widehat{\phi}(\xi)d\xi=\gamma\Big(s-\frac n {\deg p}\Big)  \sum_{i=1}^l c_{ij}(s) \int_{V_i} |p^{-s}(m)|^{-s} \phi(m) \d m,
\end{equation}
where $\gamma(s)$ is given by a product of $\Gamma$-functions, and the $c_{ij}$ are entire functions. Originally, the theory of prehomogeneous vector spaces developed from an attempt to construct Dirichlet series satisfying  functional equations in a systematic way. Since every smooth affine algebraic variety $M$ which carries the action of a reductive algebraic group $G$, and posesses an open orbit, is isomorphic to a homogeneous vector bundle $G \times _H V$, where $H$ is  a reductive subgroup of $G$, and $V$ is a locally homogeneous $H$-vector bundle, the investigation of such $G$-varieties can be reduced to the study of locally transitive linear actions of connected reductive groups.

 In classical representation theory, for every admissible irreducible representation $(\sigma,H)$ of $G$ on some Hilbert space $H$, the continuous linear operators $\sigma(f)$, defined by $\int f(g) \sigma(g) \d_{G}(g)$, where $f$ is some test function on $G$,  and $d_G$ denotes Haar measure on $G$, are of trace class, and the global character of $\sigma$ is  defined in terms of this trace as a distribution on $G$.
In this paper, we will consider the left regular representation $\pi$ of the group $G_\R$ on the Banach space $\Cvan(V_\R)$ of continuous functions on $V_\R$ vanishing at infinity, and our main goal  will consist in giving a description of the corresponding operators $ \pi(f)$, where $f$ is a rapidly decreasing function on $G$. Since we will work in a $\Cinft$ framework, avoiding Hilbert space theory entirely, we will rely on the theory of pseudodifferential operators in order to formulate our results.  It turns out that,  for arbitrary prehomogeneous vector spaces, the restriction of $\pi(f)$ to $\CT(V_\R-S_\R)$ is a pseudodifferential operator with smooth kernel. In case that $G_\C$ is reductive, and both $S_\C$ and $S_\R$ are irreducible hypersurfaces, each of the components $V_i$ is a $\Cinft$ manifold with boundary $\gd V_i$, the latter being smooth outside its intersection with the set of non-regular points $S_\R^{\rm{sing}}$ of $S_\R$.
 As we will  show, the restriction of $\pi(f)$ to $\overline V_i-S_\R^{\rm{sing}}$ is a totally characteristic pseudodifferential operator lying in the class $\L^{-\infty}_b(\overline V_i-S_\R^{\rm{sing}})$. These operators are locally given  as oscillatory integrals 
 \begin{equation*}
   Au(x)=\int e ^{i(x-y)\cdot \xi}  a(x,\xi) u(y) \d y \dbar \xi,
 \end{equation*}
where $ a(x,\xi)=\tilde a(x,x_1\xi_1,\xi')$, and $\tilde a\in \Symsl(Z\times \R^n)$ is a lacunary symbol. Here $x=(x_1,x')$ are  the standard coordinates in $Z=\overline{\R^+}\times \R^{n-1}$. The kernels of such operators are no longer smooth, and become singulare along $S_\R\times S_\R$. In the case of a smooth operator acting on a $\Cinft$ manifold, a trace can be defined in a natural way by restricting its kernel to the diagonal, which yields a density on the underlying manifold. If the manifold is compact, it can be integrated, and the so defined trace coincides with the usual $\L^2$-trace. In our situation, the restriction of the Schwartz kernel of $\pi(f)$ to the diagonal defines a distribution on $V_\R-S_\R$ given by  expressions of the form \eqref{intr1} for some critical exponent $s$, and its extension to $V_\R-S_\R^{\rm{sing}}$,   and the complex $s$-plane,   satisfies functional equations  according to \eqref{intr2}. A trace of $\pi(f)$ can then be defined by subtracting the singular contributions of the poles of the meromorphic extension.

As an application, we consider the holomorphic semigroup $S_t$ of bounded linear operators on $\Cvan(V_\R)$ generated by a strongly elliptic differential operator associated with the representation $\pi$. The latter is a differential operator of Euler type, and the corresponding semigroup can be characterized by a convolution semigroup of complex measures, which are absolutely continuous with respect to Haar measure. Denoting the corresponding Radon-Nikodym derivative by $K_t(g)\in \L^1(G)$, one has $S_t=\pi(K_t)$, and since $K_t(g)$ is analytic in $t$ and $g$, as well as rapidly decreasing, we can apply the above considerations. In particular, we get explicit expressions for the Schwartz kernel of the operators $S_t:\CT(V_\R) \rightarrow \D'(V_\R)$, and their  restrictions to the diagonal.

\section{Linear algebraic groups and prehomogeneous vector spaces}

Let $G_\C\subset \GL(m,\C)$ be a connected, linear algebraic group defined over $\R$, and $G_\R=G_\C\cap \GL(m,\R)$ the group of $\R$-rational points of $G_\C$. We will assume that $G_\C$ is symmetric, so that $G_\R$ becomes a real reductive algebraic group in the sense of \cite{wallach}, page 42. Let $\g$ be the Lie algebra of $G_\R$,     $K$ a maximal compact subgroup of $G_\R$ with Lie algebra $\k$, and  
\begin{displaymath}
\g=\k\oplus\p
\end{displaymath}
the corresponding Cartan decomposition of $\g$.  The mapping $(k,X)\mapsto k\exp X$ represents an analytic diffeomorphism of $K\times \p$ onto $G_\R$, while $\exp:\k\rightarrow K$ is a surjection. Let $G$ be the connected component of $G_\R$ containing the unit element $e$. Clearly, $G$ is a Lie group, and its Lie algebra coincides with $\g$. 
Choose a left invariant Riemannian metric on $G_\R$, and denote the distance of two points $g,h \in G$  by $\d(g,h)$. We set $|g|=d(g,e)$. Then
\begin{displaymath}
|g|=|g^{-1}|, \qquad |e|=0, \qquad |gh|\leq |g|+|h|, \qquad g,h \in G.
\end{displaymath}
In the following, we will say that a function $f$ on $G$ is \emph{at most of exponential growth}, if there exists a $\kappa>0$ such that $|f(g)| \leq C e^{\kappa|g|}$ for some constant $C>0$.
If $g=k\exp X$ is the Cartan decomposition of an arbitrary element in $G$, one computes             
\begin{equation*}
|g| \leq |k|+|\exp X|\leq C_K +|\exp X|, \qquad |\exp X|=|k^{-1}g|\leq
 |k|+|g|\leq C_K+|g|,
\end{equation*}
where $C_K=\max_{k\in K} |k|<\infty$. Putting $C'=e^{C_K}\geq 1$, we obtain 
\begin{equation}
\label{F}
\frac 1 {C'} e^{|\exp X|}\leq e^{|g|}\leq C'e^{|\exp X|}.
\end{equation}
Let $d$ be the dimension of $G_\R$, $X_1,\dots,X_l$ be a basis of  $\k$ and $X_{l+1},\dots,X_d$ a basis
of $\p$. Since $\exp:\p\simeq P$ is an analytic diffeomorphism, there exists a $C''\geq 1$, such that
\begin{equation}
\label{G}
\frac 1 {C''} |X|\leq|\exp X|\leq C''| X|,
\end{equation}
where $|X|=\sqrt{q_{l+1}^2+\dots+ q_d^2}$ is the length of
$X=q_{l+1}X_{l+1}+\dots+q_dX_d\in \p$.  Realizing $\exp$ as power series for matrices, relations \eqref{F} and \eqref{G} then imply that the matrix coefficients of $\exp X$ and,  consequently, of $g=k\exp X$, are at most of exponential growth.

Let $L$ respectively $ R$ be the left respectively right regular representation of $G$ on the Fr\'{e}chet space $\Cinft(G)$ of smooth, complex valued functions on $G$, equipped with the usual topology of uniform convergence on compact subsets of a function and each of its derivatives, see \cite{warner}, page  220. Let $\U$ be the universal enveloping algebra of the complexification of  $\g$, and denote the representations of $\U$ on the space  $\Cinft(G)_\infty$ of differentiable elements  by $dL$, respectively $dR$. Let $d_{G}$ denote left invariant Haar measure on $G$. We introduce now the space of rapidly decreasing $\Cinft$ functions on $G$.

\begin{definition}
\label{def:A}
The \emph{space of rapidly decreasing $\Cinft$ functions on $G$}, in the following denoted by  $\S(G)$,  is given by all functions $f\in\Cinft(G)_\infty$ satisfying the following conditions:

\noindent
(i) For every $\kappa\geq 0$, and $X\in \U$, there exists a constant $C$  such that 
\begin{displaymath}
|dL(X) f(g)|\leq C e^{-\kappa |g|};
\end{displaymath}

\noindent
(ii) for every   $\kappa\geq 0$, and $X\in \U$, one has $dL(X)f\in
\L^1(G,e^{\kappa |g|}d_G)$. 

\noindent

\end{definition}

Note that $\S(G)$ is $\pi(G)$- and $d\pi(\U)$-invariant.

\begin{remark}
\label{rem:A}
Let  $f\in \S(G)$. By the relation $g\e{X}g^{-1}=\e{\Ad(g)X}$, $X \in \g$, and with respect to the basis of $\g$ introduced above, one calculates
\begin{align*}
dR(X_i)f(g)&=\lim_{h \to 0} h^{-1} [f(g\e{hX_i}g^{-1} g)-f(g)]=-dL(\Ad(g) X_i)f(g)\\&=-\sum \limits _{j=1}^d A_{ij}(g)dL(X_j)f(g),
\end{align*}
where the matrix coefficients $A_{ij}$ of the adjoint representation, and their derivatives, are at most of exponential growth. Thus, for  arbitrary $X\in \U$ and $\kappa\geq 0$, $dR(X)f$ satisfies the conditions $(i)$ and $(ii)$ of the preceeding definition.
\end{remark}
In the sequel, the following partial integration formulas will be needed. Let  $\Delta(g)=|\det \Ad(g)|$ be  the modular function of $G$. 
\begin{proposition}
\label{prop:A}
Let $f_1\in\S(G)$, and assume that $f_2 \in \Cinft(G)_\infty$, together with all its derivatives, is at most of exponential growth.  Then, for arbitrary multiindices $\gamma$, 
\begin{equation}
\label{A}
\int_{G}f_1(g) dL(X^\gamma) f_2(g) d_{G}(g)=(-1)^{|\gamma|} \int _{G}
dL(X^{\tilde \gamma}) f_1(g) f_2(g) d_{G}(g),
\end{equation}
\begin{equation}
\label{B}
\int_{G}f_1(g) dR(X^\gamma) f_2(g)
d_{G}(g)=\sum\limits_{\alpha_1+\alpha_2=\gamma} \iota ^{\alpha_1} (-1)^{|\alpha_2|} \int _{G}
dR(X^{\tilde \alpha_2}) f_1(g) f_2(g) d_{G}(g),
\end{equation} 
where   $X^\gamma=X^{\gamma_1}_{i_1}\dots X^{\gamma_r}_{i_r}$,   
$X^{\tilde \gamma}=X^{\gamma_r}_{i_r}\dots X^{\gamma_1}_{i_1}$, and   $\iota_k=dR(X_k)\Delta(e)=dL(-X_k)\Delta(e)$.   
\end{proposition}
\begin{proof}

First, note that $f_1(g) f_2(\e{-hX_i}g)$ is  a differentiable function with respect to $h$, and integrable on $G$ for all $h \in \R$, since $|f_1(g) f_2(\e{-hX_i}g)|\leq C |f_1(g)| e^{\kappa (|\e{hX_i}|+|g|)}$ for some $C,\kappa >0$, and $f_1(g)\in \L^1(G, e^{\kappa|g|}d_G)$. Furthermore,
\begin{displaymath}
\frac d{dh} f(\e{-hX_i}g)_{|h=h_0}=\frac d{dh}
f(\e{-hX_i}\e{-h_0X_i}g)_{|h=0}=\pi(\e{h_0X_i})dL(X_i)f(g),
\end{displaymath}
so that, for  $h_0\in \R$,
\begin{align*}
\modulus{f_1(g) \frac d {dh} f_2(\e{-hX_i}g)_{|h=h_0}}=|f_1(g)
\pi(\e{h_0X_i})dL(X_i)f_2(g)|\leq C' |f_1(g)|e^{\kappa' (|\e{hX_i}|+|g|)},
\end{align*}
where $C',\kappa' >0$. Thus, for $h_0 \in (-\epsilon,\epsilon)$, $f_1(g) \frac d {dh} f_2(\e{-hX_i}g)_{|h=h_0}$ can be estimated from above by $C_\epsilon |f_1(g)| e^{\kappa'}$, $C_\epsilon>0$ being some constant. By the bounded convergence theorem of Lebesgue, and the left invariance of Haar measure, we therefore obtain
\begin{gather*}
\int_{G} f_1(g) dL(X_i) f_2(g)d_{G}(g)=\int_{G}\lim _{h \to 0}
h^{-1}[f_1(g) (f_2(\e{-hX_i}g)-f_2(g))] d_G(g)\\
=\lim_{h \to 0} h^{-1} \int _{G} f_1(g) [f_2(\e{-hX_i}g)-f_2(g)] d_{G}(g)=\lim_{h \to 0} h^{-1} \int _{G}  [f_1(\e{hX_i}g)-f_1(g)]f_2(g) d_{G}(g),
\end{gather*}
 compare \cite{elstrodt}, pages 146 and 364. Since $f_1(\e{hX_i}g)f_2(g)$ and $f_2(g) \frac d {dh}
f_1(\e{hX_i}g)_{|h=h_0}$, $h_0 \in (-\epsilon,\epsilon)$ can be  estimated in a similar way, a repeated application of Lebesgue's theorem finally yields  \eqref{A}
for   $|\gamma|=1$. Similarly, by taking into account Remark \ref{rem:A}, we deduce
\begin{gather*}
\int_{G} f_1(g) dR(X_i) f_2(g)d_{G}(g)=
\lim_{h \to 0} h^{-1} \int _{G} f_1(g) [f_2(g\e{hX_i})-f_2(g)] d_{G}(g)\\
=\lim_{h \to 0} h^{-1} \int _{G}  [f_1(g\e{-hX_i})\Delta(\e{-hX_i})^{-1}-f_1(g)]f_2(g) d_{G}(g),
\end{gather*}
obtaining \eqref{B} for $|\gamma|=1$. The general formulas then follow by induction. 
\end{proof}
 In what follows, we will denote the coordinate functions in $M_n(\R)\simeq \R^{n^2}$ by  $g_{ij}$.
 With  the identification $\g\simeq \R^d$, and with respect to a basis $X_1, \dots, X_d$ of $\g$, the canonical coordinates of second type of a point $g\in G_\R$ are given by 
\begin{equation}
\label{C}
\Phi_g:g U_e\ni g \e{\zeta_1a_1}\dots \e{\zeta_da_d} \mapsto (\zeta_1,\dots, \zeta_d) \in W_0,
\end{equation}
where $W_0$ denotes a sufficiently small neighbourhood of $0$ in $\R^d$, and $U_e=\exp(W_0)$. We will write for $\Phi_e$ simply  $\Phi$. As a real analytic submanifold of $\GL(n,\R)$, $G_\R$ is, in particular, embedded in  $M_n(\R)$, so that the matrix
\begin{equation*}
\frac {D(g_{11}\circ \Phi_g^{-1},\dots,g_{nn}\circ \Phi_g^{-1})}{D(\zeta_1,\dots, \zeta_d)}=\left ( \frac d {d\zeta_k} (g_{ij}\circ \Phi_g^{-1})(\zeta)\right )_{ij,k}
\end{equation*}
has rang $d$ for all $g \in G_\R$ and $\zeta \in W_0$. For each tuple ${i_1j_1,\dots,i_dj_d}$ of indices we then define the sets 
\begin{equation}
\label{D}
  G_\R^{i_1j_1,\dots,i_dj_d}=\mklm{g \in G_\R: \det \frac {D(g_{i_1j_1}\circ \Phi_g^{-1},\dots,g_{i_dj_d}\circ \Phi_g^{-1})}{D(\zeta_1,\dots, \zeta_d)}\not=0},
\end{equation}
in this way getting a finite covering of $G_\R$ by open, although  not necessarily connected, sets.

Let $V_\C$ be a $n$-dimensional vector space over $\C$ with $\R$-structure $V_\R$, $\rho:G_\C\rightarrow \GL(V_\C)$ a $\R$-rational representation of $G_\C$ on $V_\C$, and suppose that $\rho$ has a dense $G_\C$-orbit in the Zariski topology. Then $(G_\C,\rho,V_\C)$ is called a \emph{prehomogeneous vector space over $\R$}. The complement set of the unique dense orbit, the \emph{singular set}, is a Zariski closed set, and will be denoted by $S_\C$.  Since $\rho(G_R)\subset \rho(G)_\R$, the restriction of $\rho$ to $G_\R$ induces a regular $G_\R$-action on $V_\R$.  Assume now that $G_\C$ is reductive, and that the singular set $S_\C$ is an irreducible hypersurface. Then, there is an irreducible, homogeneous polynomial $p$ such that $S_\C=\mklm{m \in V_\C: p(m)=0}$.
It  is called a \emph{relative invariant}, and there exists a rational character $\chi:G_\C \rightarrow \GL(1,\C)$, such that
\begin{equation}
 \label{E}
  p(\rho(g)m))=\chi(g) p(m), \qquad m \in V_\C-S_\C, \quad  g \in G_\C.
\end{equation}
In this case, $(G_\C,\rho,V_\C)$ is  a regular prehomogeneous vector space, and one has $\deg p | 2n$, $\det \rho(g)^2=\chi(g) ^{2n/\deg p}$, where $n=\dim_\C V_\C$. By multiplying $p$ with a scalar, we can assume that $p(V_\R)\subset \R$, and $\chi(G_\R) \subset \R^\ast$. Let $S_\R=S_\C \cap V_\R=\mklm{m \in V_\R: p(m)=0}$.  As before, denote by $G$ the connected component of $G_\R$ containing the unit element. Then, by a theorem of Whitney, $V_\R-S_\R$ decomposes into a finite number of connected components, each of them being a $G$-orbit, and we write $V_\R-S_\R=V_1 \cup \dots \cup V_l$. In what follows, we will identify $V_\C$ with $\C^n$, and $V_\R$ with $\R^n$, by choosing a basis in $V_\R$, and assume that $G \subset \GL(n,\C)$ without loss of generality. Instead of $\rho(g)m$ we will then simply write $gm$. In particular, $S_\C$ and $S_\R$ become irreducible affine algebraic varieties in $\C^n$, respectively $\R^n$. Their sets of regular points, $S_\C^{\rm{reg}}$, respectively  $S_\R^{\rm{reg}}$ can be provided  with differentiable structures, the underlying topology being the induced one. Let $(G_\C,\rho^\ast,V_\C^\ast)$ denote the dual prehomogeneous vector space of $(G_\C,\rho,V_\C)$, and $p^\ast$ its relative invariant satisfying $p^\ast(\rho^\ast(g)\xi)=\chi(g)^{-1} p^\ast(\xi)$. Let $\S(V_\R)$, respectively $\S(V_\R^\ast)$, denote the Schwartz space of functions on $V_\R$, respectively $V_\R^\ast$, and let  $\hat \phi$ be  the Fourier transform of $\phi\in \S(V_\R)$.  Let $dm$ be Lebesgue measure on $V_\R$. The following result is known as the  fundamental theorem  of prehomogeneous vector spaces, and was proved by Sato in 1961 \cite{kimura}.                         
\begin{theorem}[Sato]
\label{thm:A}
Let $\phi \in \S(V_\R)$ and $\phi^\ast \in \S(V^\ast_\R)$. Then the integrals 
\begin{equation}
\label{J}
 F_j(s,\phi)= \frac 1 {\gamma(s)} \int_{V_j} |p(m)|^s \phi(m) dm, \qquad F_i^\ast(s,\phi)=\frac 1 {\gamma(s)} \int_{V_i^\ast} |p^\ast(\xi)|^s \phi^\ast(\xi) d\xi,
\end{equation}
converge for $\Re s >0$, and can be extended analytically to holomorphic functions on the whole $s$-plane, satisfying the functional equations
\begin{equation}
\label{K}
 F_i^\ast(s-n/\deg p,\hat \phi)= \gamma(-s)\sum _{j=1}^l c_{ij}(s)F_j(-s,\phi), 
\end{equation}
where $\gamma(s)$ is given by a product of $\Gamma$-functions, and  the $c_{ij}(s)$ are entire functions which do not depend on $\phi$.
The functions $F_j(s,\phi)$, $F^\ast_j(s,\phi)$ are called \emph{local $\zeta$-functions}.
\end{theorem}

\section{Review of pseudodifferential operators}
\label{sec:III}

\subsection{Generalities}

This section is devoted to the exposition of some of the basic facts about  pseudodifferential operators, needed to formulate our main results in the sequel. Our main references for the theory will be \cite{hoermanderIII} and \cite{shubin}.   Consider first an open set $U$ in $\R^n$, and let $x_1,\dots,x_n$ be the standard coordinates. For any real number $l$,  we denote by $\Sym^l(U\times \R^n)$ the class of all functions $a(x,\xi)\in \Cinft(U\times \R^n)$ such that, for any multiindices $\alpha,\beta$, and any compact set $K\subset U$, there exist constants $C_{\alpha,\beta,K}$ for which
\begin{equation}
  \label{H}
|\gd ^\alpha_\xi\gd ^\beta_x a(x,\xi)| \leq C_{\alpha,\beta,K} \eklm{\xi}^{l-|\alpha|}, \qquad x \in K,  \quad \xi \in \R^n,
\end{equation}
where $\eklm{\xi}$ stands for $(1+|\xi|^2)^{1/2}$, and $|\alpha|=\alpha_1,\dots,\alpha_n$. We further put $\Sym^{-\infty}(U\times \R^n)=\bigcap _{l \in \R} \Sym^l(U\times \R^n)$. Note that, in general, the constants $C_{\alpha,\beta,K}$ also depend on $a(x,\xi)$. For any such $a(x,\xi)$ one then defines the continuous linear operator
\begin{displaymath}
  A:\CT(U) \longrightarrow \Cinft(U)
\end{displaymath}
by the formula
\begin{equation}
  \label{I}
Au(x)=\int e^{ix \cdot \xi} a(x,\xi) \hat u(\xi) \dbar \xi,
\end{equation}
where $\hat u$ denotes the Fourier transform of $u$, and $\dbar \xi=(2\pi)^{-n} \d \xi$. An operator $A$ of this form is called a \emph{pseudodifferential operator}, and we denote the class of all such operators for which $a(x,\xi) \in \Sym^l(U\times \R^n)$  by $\L^l(U)$. The set $\L^{-\infty}(U)=\bigcap_{l\in \R} \L^l (U)$ consists of all operators with smooth kernel. Such operators are called \emph{smooth operators}. 
By inserting  in \eqref{I} the definition of $\hat u$, we obtain for $Au$ the expression
\begin{equation}
  \label{II}
Au(x)=\int \int  e^{i(x-y) \cdot \xi} a(x,\xi)  u(y) \d y \,\dbar \xi, 
\end{equation}
which exists as an oscillatory integral, see \cite{shubin}. The Schwartz kernel $K_A\in \D '(U\times U)$ of $A$ is therefore given by the oscillatory integral
\begin{equation}
  \label{III}
K_A(x,y)=\int e^{i(x-y)\cdot \xi} a(x,\xi) \,\dbar \xi,
\end{equation}
and is a smooth function outside the diagonal in $U\times U$.

Consider next an $n$-dimensional $\Cinft$ manifold $X$, and let  $(\kappa_\gamma, U^\gamma)$ be an atlas for $X$. Then a linear operator
\begin{equation}
\label{IIIa}
  A:\CT(X) \longrightarrow \Cinft(X)
\end{equation}
is called a \emph{pseudodifferential operator on $X$ of order $l$} if, for every chart diffeomorphism $\kappa_\gamma:U^\gamma \rightarrow \tilde U^\gamma= \kappa_\gamma(U^\gamma)$,  the operator $A^{\gamma} u = [A_{|U^\gamma} ( u\circ \kappa_{\gamma})] \circ \kappa_{\gamma}^{-1}$ given by the diagramm
\begin{diagram}
\CT(U^{\gamma})           & \rTo^{\qquad A_{|U^\gamma}}             & & \Cinft(U^\gamma)               \\
\kappa_{\gamma}^\ast \, \uTo &                                  & & \uTo(2,1)\, \kappa_\gamma^\ast \\
\CT(\tilde U^{\gamma})    & \rTo^{\quad \quad A^{\gamma}} & & \Cinft(\tilde U^\gamma)        \\
\end{diagram}
is a pseudodifferential operator on $\tilde U^\gamma$ of order $l$, and we write $A \in \L^l(X)$. Note that, since the $U^\gamma$ are not necessarily connected, we can choose them in such a way that $X\times X$ is covered by the open sets $U^\gamma \times U^\gamma$. 
Now, in general, if $X$ and $Y$ are two smooth manifolds, and 
\begin{equation*}
  A: \CT(X) \longrightarrow \Cinft(Y) \subset \D'(Y)
\end{equation*}
is a continuous linear operator, where $\D'(Y)=(\CT(Y,\Omega))'$ and $\Omega=|\Lambda^n(Y)|$ is the density bundle on $Y$, its Schwartz kernel is given by the distribution section $K_A \in \D'(Y \times X, \1 \boxtimes \Omega_X)$, where  $\D'(Y\times X ,1 \boxtimes \Omega_X) = (\CT(Y \times X, (\1 \boxtimes \Omega_X)^\ast \otimes \Omega_{Y\times X}))'$. Observe that  $\CT (Y, \Omega_Y) \otimes \Cinft (X) \simeq \Cinft ( Y \times X, (\1 \boxtimes \Omega_X)^\ast \otimes \Omega_{Y\times X})$. In case that  $X=Y$   and $A\in \L^l(X)$, $A$ is given locally by the operators $A^{\gamma}$, which can be written in the form 
  \begin{equation*}
    A^{\gamma}u(x) = \int e^{i(x-y) \cdot \xi} a ^{\gamma}(x,\xi) u(y) \d y \dbar \xi,
  \end{equation*}
where $u \in \CT(\tilde U^{\gamma})$, $x \in \tilde U^\gamma$, and $a^{\gamma}(x,y,\xi) \in {\rm{S}}^l(\tilde U^\gamma, \R^n)$. The kernel of $A$ is then determined by the kernels $K_{A^{\gamma}} \in \D'(\tilde U^\gamma \times \tilde U ^{\gamma})$. For $l < -\dim X$, they are continuous and given by absolutely convergent integrals. 
In this case, their restrictions to the respective diagonals in $\tilde U^\gamma \times \tilde U^\gamma$ define continuous functions
\begin{equation*}
  k^\gamma(m)= K_{A^{\gamma}} (\kappa_\gamma (m),\kappa _\gamma (m)), \qquad m \in U^\gamma,
\end{equation*}
which, for $m \in  U^{\gamma_1} \cap U^{\gamma_2}$, satisfy the relations $ k^{\gamma_2}(m) =| \det (\kappa_{\gamma_1} \circ \kappa ^{-1}_{\gamma_2})' | \circ \kappa_{\gamma_2}(m) k^{\gamma_1}(m)$, thus defining a density $k \in \Cinft (X,\Omega)$ on $\Delta_{X\times X} \simeq X$.   


\subsection{Totally characteristic pseudodifferential operators} We introduce now a special class of pseudodifferential operators associated in a natural way to a $\Cinft$ manifold $X$ with boundary $\gd X$.  Our main reference will be \cite{melrose} in this case. Let $\Cinft(X)$ be the space of functions on $X$ which are $\Cinft$ up to the boundary, and $\Cinftv(X)$ the subspace of functions vanishing to all orders on $\gd X$. The standard spaces of distributions over $X$ are 
\begin{equation*}
  \D'(X)= (\CTv(X,\Omega))', \qquad \dot D'(X) =(\CT(X,\Omega))', 
\end{equation*}
the first being the space of extendible distributions, whereas the second is the space of distributions supported by $X$. Consider now the translated partial Fourier transform of a symbol $a(x,\xi) \in \Sym^l(\R^n\times \R^n)$,
\begin{equation}
  Ma(x,\xi';t)=\int e^{i(1-t)\xi_1} a(x,\xi_1,\xi') d\xi_1,
\end{equation}
which is $\Cinft$ away from $t=1$. Here $\xi=(\xi_1,\xi')$. We will say that $a(x,\xi)$ is \emph{lacunary} if it satisfies the  \emph{lacunary condition}
\begin{equation}
\label{V}
  Ma(x,\xi';t)=0 \qquad \text{ for } t<0.
\end{equation}
The subspace of lacunary symbols will be denoted by $\Sym^l_{la}(\R^n\times \R^n)$. Let $Z=\overline{\R^+}  \times \R^{n-1}$ be the standard manifold with boundary  with the natural coordinates $x=(x_1,x')$. In order to define on $Z$ operators written formally as oscillatory integrals
\begin{equation}
\label{VI}
  Au(x)=  \int e^{i(x-y)\xi }  a (x,\xi) u(y) \, \d y \dbar \xi, 
\end{equation}
where $a (x,\xi)=\tilde a(x_1,x',x_1\xi_1, \xi')$ is a more general amplitude and $\tilde a(x,\xi)$ is lacunary, one rewrites the formal adjoint of $A$ by making a singular coordinate change. Thus, for $u \in \CT(Z)$, one considers 
\begin{equation*}
  A^\ast u(y) =\int e^{i(y-x) \xi} \overline{a}(x,\xi) u(x) \, \d x \dbar \xi.
\end{equation*}
Putting $\lambda=x_1\xi_1$, $s=x_1/y_1$ one obtains
\begin{equation}
\label{VII}
  A^\ast u(y)=(2\pi)^{-n}\int e^{i(1/s-1,y'-x')\cdot (\lambda,\xi')}\overline {\tilde a}(y_1 s, x',\lambda,\xi') u(y_1 s,x')d\lambda \frac{ds}s dx' d\xi'.
\end{equation}
For $\tilde a \in \Sym_{la}^{-\infty}(Z\times \R^n) $, the succesive integrals in \eqref{VII} converge absolutely and uniformly, thus defining  a continuous bilinear form
  \begin{equation*}
    \Sym^{-\infty}_{la}(Z\times \R^n) \times \CT(Z) \longrightarrow \Cinft(Z),
  \end{equation*}
which by  \cite{melrose}, Propositions 3.6 and 3.9. extends to a separately continuous form
\begin{equation*}
  \Sym^\infty_{la}(Z\times \R^n) \times \CT(Z) \longrightarrow \Cinft(Z).
\end{equation*}
If $\tilde  a \in \Sym^\infty_{la}(Z \times \R^n)$, one then defines the operator 
\begin{equation}
\label{VIII}
  A:\dot \E'(Z) \longrightarrow \dot \D'(Z), 
\end{equation}
written formally as \eqref{VI}, as the adjoint of $A^\ast$. In this way, the oscillatory integral \eqref{VI} is identified with a separately continuous bilinear mapping
\begin{equation*}
  \Sym^\infty_{la}(Z\times \R^n) \times \dot \E'(Z) \longrightarrow \dot \D'(Z).
\end{equation*}
 The space $\L^l_b(Z)$ of \emph{totally characteristic pseudodifferential operators on $Z$ of order $l$} consists of those continuous linear maps \eqref{VIII} such that $v Au$ is of the form \eqref{VI} with $\tilde a(x,\xi)\in \Sym^l_{la}(Z\times \R^n)$ whenever $u,v \in \CT(Z)$. Similarly, 
a continuous linear map \eqref{IIIa} on a smooth manifold $X$ with boundary $\gd X$ is an element of the space $\L^l_{b}(X)$ of totally characteristic pseudodifferential operators on $X$  if, for a given atlas $(\kappa_\gamma,U^\gamma)$, the operators $A^\gamma u=[A_{|U^\gamma} (u \circ \kappa_\gamma)] \circ \kappa^{-1}_\gamma$ are elements of $\L^l_b(Z)$,  
where $U^\gamma$ are coordinate patches isomorphic to subsets in $Z$.

\section{A structure theorem for prehomogeneous vector spaces}

\label{sec:IV}

Let $(G_\C,\rho,V_\C)$ be a prehomogeneous vector space defined over $\R$, and identify $V_\C$ with $\C^n$, and $V_\R$ with $\R^n$, by choosing a basis in $V_\R$. Consider the Banach space $\Cvan(V_\R)$ of continuous, complex valued functions on $V_\R$ vanishing at infinity, equipped with the supremum norm. Let  $(\pi,\Cvan(V_\R))$ be the corresponding, continuous,  left regular representation of $G_\R$, which  was already introduced in \cite{ramacherI}.  The representation of $\U$ on the space of differentiable vectors $\Cvan(V_\R)_\infty$ will be denoted by $d\pi$. We will also consider the left regular representation  of $G_\R$ on $\Cinft(V_\R)$ which, equipped with the topology of uniform convergence on compact subsets, becomes a Fr\'{e}chet space. This representation will be denoted by $\pi$ as well. As before, we write $G$ for the connected component of $G_\R$ containing the unit element, and $d_G$ for Haar measure on $G$.  Let  $f\in \S(G)$, and denote by $\pi(f)$ the   continuous linear operator $\int _{G} f(g)  
\pi(g) \d_{G}(g)$ in $\Cvan(V_\R)$. Its restriction to $\CT(V_\R)$ induces a continuous linear operator            
\begin{equation}
\label{-1}
\pi(f):\CT(V_\R) \longrightarrow \Cvan(V_\R) \subset \D'(V_\R).
\end{equation}
To see this, let  $u_j\in\CT(V_\R)$ be a sequence which converges to zero in $\CT(V_\R)$, that is, assume that there exists a compact set $K\subset V_\R$ such that $\supp u_j \subset K$, and $\sup |\gd ^\alpha u_j|\to 0$ for all   $j$ and arbitrary multiindices $\alpha$. Then, for  $v\in \CT(V_\R)$, one has
\begin{displaymath}
|\eklm{\pi(f) u_j,v}|=\modulus{\int_{V_\R} (\pi(f) u_j)(m) v(m) dm} \leq \norm{\pi(f)u_j} \cdot \norm{v}_{L^1} \to 0,
\end{displaymath}
since $\norm{\pi(f) u_j}\leq\norm{\pi(f)}\norm{u_j}$, and $\norm{u_j}\to 0$ by assumption. As a consequence, and according to Schwartz, there exists          a distribution $\K_f \in \D'(V_\R\times V_\R)$ such that $\eklm{\pi(f)u,v}=\K_f(v\otimes u)$ for all   $u,v \in
\CT(V_\R)$. The  properties of the Schwartz kernel $\K_f$ will depend    on the analytic properties of $f$, and the orbit structure of the underlying $G$-action, and our main effort in the following sections will be directed towards the elucidation of the structure of $\K_f$. Denote the coordinates of a point in $V_\R$, respectively $V_\R^\ast$, by $m_1, \dots, m_n$, respectively $\xi_1,\dots,\xi_n$, and put       
$$ \phi_\xi(m)=e^{i m\cdot \xi}=e^{i\sum_{j=1}^n m_j
\xi_j},\qquad m\in V_\R, \, \xi \in V_\R^\ast.$$
Clearly, $\phi_\xi(m)\in \Cinft(V_\R\times V_\R^\ast)$. Since $f\in
\L^1(G,d_{G})$, we can define the function
\begin{equation}
\label{0}
\hat f(m,\xi)=\int _{G} f(g) (\pi(g)\phi_\xi)(m) d_{G}(g)=\int _{G} f(g) e^{i (g^{-1}m)\cdot \xi} d_{G}( g), \quad
m\in V_\R,\,\xi \in V_\R^\ast.
\end{equation}
$\hat f$ is continuous in $m$ and $\xi$, and  one verifies $|\hat f(m,\xi)|\leq \norm{f}_{\L^1}$. Since  $\gd
^\alpha_\xi\gd^\beta_m\pi(g) \phi_\xi(m)$ is,  together with all its derivatives,   at most of exponential growth in $g$,
 we can interchange  the order of integration and differentiation  in \eqref{0}, yielding $\hat f (m,\xi)\in
\Cinft(V_\R\times \R_\xi^n)$. The main issue of this section will consist in proving the following structure theorem.

\begin{theorem}
\label{thm:1}
Let $(G_\C, \rho, V_\C)$ be a prehomogeneous vector space over $\R$, and $(\pi,\Cvan(V_R))$ the left regular representation of $G_\R$ on $V_\R$. Then, for $f\in \S(G)$, $\pi(f):\CT(V_\R)\rightarrow \D'(V_\R- S_\R)$ is a Fourier integral operator given by         
\begin{equation}
\label{1}
(\pi(f) \phi) (m)=\int \int e^{i(m-m')\cdot\xi} a_f(m,\xi) \phi(m') \d
m' \,\dbar \xi,
\end{equation}
with symbol
\begin{equation*}
a_f(m,\xi)=e^{-im\cdot \xi} \hat f(m,\xi) \in \Syms((V_\R- S_\R)\times V_\R^\ast).
\end{equation*}
On  $V_\R- S_\R$, $\pi(f)$ is a pseudodifferential operator of class $\L^{-\infty}(V_\R
- S_\R)$.
\end{theorem}

For the proof, we will need some lemmas. First, note that on $V_\R-S_\R$, one can express the canonical vector fields $\gd_{m_i}$ of  $V_\R$ locally as linear combinations in the fundamental vector fields of the underlying $G_\R$-action.

\begin{lemma}
\label{lem:1}
Let $\tilde m \in V_\R- S_\R$. Then there 
exists a neighbourhood $U\subset V_\R- S_\R$ of  $\tilde m$, and rational functions  $\Theta^k_i(m)$ in $m$ on  $U$, such that, for arbitrary $\phi \in \Cinft(V_\R)$, 
\begin{equation}
(\gd_{m_i} \phi)(m)=\sum_{k=1}^n \Theta^k_i(m)
(d\pi(X_{j_k})\phi)(m), \qquad m \in U,
\end{equation}
for some indices $j_1,\dots,j_n$.
\end{lemma}
\begin{proof}
For $m \in V_\R$, one has
\begin{equation}
\label{1a}
(d\pi(X_j)\phi)(m)=\frac d {dh} \phi(e^{-hX_j}m)_{\big\vert_{h=0}}=\sum_{i=1}^n
 \frac d{dh} m_i(e^{-hX_j}m)_{\big\vert_{h=0}}(\gd_{m_i} \phi)(m).
\end{equation}
Because of the local transitivity of the $G_\R$-action on $V_\R-
S_\R$, the $n\times d$-matrix
\begin{displaymath}
\Big (\frac d{dh} m_i(e^{-hX_j}\tilde m)_{\big\vert_{h=0}}\Big )_{i,j}=\Big ((d\pi(X_j)
m_i)(\tilde m) \Big )_{i,j}
\end{displaymath}
has maximal rang $n$ for $\tilde m\in V_\R- S_\R$, where, by assumption, $d\geq n$. Therefore, there exists a neighbourhood $U$ of 
$\tilde m$, and indices  $j_1,\dots,j_n$, such that for all $m \in U$
\begin{equation*}
\det \left (\frac d{dh} m_i(e^{-hX_{j_k}}m)_{\big\vert_{h=0}}\right )_{i,k} \not=0.
\end{equation*}
Denoting the matrix coefficients of the corresponding inverse matrix by $\Theta^k_i(m)$, the assertion follows with \eqref{1a}.
\end{proof}
Consider the covering of $G_\R$ by the sets $G_\R^{i_1j_1,\dots,i_dj_d}$ introduced in \eqref{D}, and put $G^{i_1j_1,\dots i_dj_d}=G_\R^{i_1j_1,\dots,i_dj_d}\cap G$.
\begin{proposition}
\label{prop:1}
Fix $\tilde m \in V_\R- S_\R$. Then there exists a neighbourhood $U\subset V_\R- S_\R$ of  $\tilde m$, such that, for arbitrary multiindices $\beta$,  $\phi \in \Cinft(V_\R)$, and $g \in G$,
\begin{equation}
\label{2}
(\gd^\beta_m\pi(g)\phi)(m)=\sum\limits_{|\gamma|\leq |\beta|} b^\gamma_\beta(m,g) \pi(g)d\pi(X^\gamma)\phi(m),\qquad m \in U.
\end{equation}
The coefficient functions $b^\gamma_\beta(m,g)\in \Cinft(U\times G)$  are rational expressions in the coordinates of  $m$ and satisfy the estimates
\begin{equation}
\label{3}
|\gd _m^\alpha b^\gamma_\beta(m,g)|\leq C_U e^{\kappa |g|}
\end{equation}
for some appropriate $C_U,\kappa>0$.
\end{proposition}
\begin{proof}
  By the previous lemma, there exists a neighbourhood $U$ of  $\tilde m$ such that, for $m \in U$ and all $g \in G$,
\begin{displaymath}
\gd_{m_i} \pi(g) \phi(m)=\sum \limits_{k=1}^d \Theta^k_i(m) d\pi(X_k) \pi(g) \phi(m).
\end{displaymath}
Now,  
\begin{displaymath}
d\pi(X_k)\pi(g) \phi(m)=\sum \limits_{j=1}^d  \frac d {dh} s^k_j(h,g)_{\big\vert_{h=0}}d\pi(X_j)\phi(g^{-1}m),
\end{displaymath}
where the functions $s^k_j$ are given by the equations
$\e{hX_k}g=\Phi_g^{-1} (s_1^k(h,g), \dots, s_d^k(h,g))=g\e{s^k_1(h,g)X_1}\dots \e{s^k_d(h,g)X_d}$, and are real analytic in $h$ and $g$, compare \cite{ramacherI}. Differentiating these equations with respect to $h$ one obtains 
\begin{equation*}
g_{ij}(X_k g)=\sum \limits_{l=1}^d \frac \gd{\gd \zeta_l} (g_{ij} \circ \Phi^{-1}_g)(\zeta)_{|\zeta=0}  \frac d {dh} s^k_l(h,g)_{\big\vert_{h=0}}.
\end{equation*}
 Let $\chi_{i_1j_1,\dots,i_dj_d}(g) \in \Cinft(G)$ be a finite partition of unity subordinate to the covering ${G^{i_1j_1,\dots,i_dj_d}}$ of $G$, that is, assume that $\supp \chi_{i_1j_1,\dots,i_dj_d}\subset G^{i_1j_1,\dots,i_dj_d}$,
$\chi_{i_1j_1,\dots,i_dj_d}\geq 0$, and $\sum \chi_{i_1j_1,\dots,i_dj_d}(g)=1$ for each $g\in G$.
Note that we are allowing the $\chi_{i_1j_1,\dots,i_dj_d}$ to have non-compact support.  For $g \in G^{i_1j_1,\dots, i_dj_d}$,   $D(g_{i_1j_1}\circ \Phi_g^{-1},\dots,g_{i_dj_d}\circ \Phi_g^{-1})/D(\zeta_1,\dots, \zeta_d)$ becomes invertible.  Since $\frac d {d\zeta_k} (g_{ij} \circ \Phi^{-1}_g)(\zeta)_{|\zeta =0}=g_{ij}(g X_k)$, we therefore obtain, as a consequence of the last equation, the relation
\begin{equation}
\label{4}
\left ( \begin{array}{c}  \dot s^k_1(0,g) \\ \vdots \\ \dot s^k_d(0,g)   \end{array}
\right ) =\sum _{i_1j_1,\dots, i_dj_d} \chi_{i_1j_1,\dots, i_dj_d}(g) \Big ( g_{ij}(gX_k) \Big )^{-1}_{ij,k} \left( \begin{array}{c} g_{i_1j_1}(X_k  g) \\ \vdots\\ g_{i_dj_d}(X_k g) \end{array} \right ), \quad g \in G.
\end{equation}
The coefficients of $g$ and  the functions $\chi_{i_1j_1,\dots,i_dj_d}$ are, together with their derivatives, at most of exponential growth,  so that equation \eqref{4} implies 
that the derivates of the functions $s^k_l(h,g)$ with respect to $h$ satisfy the estimates
\begin{equation}
\label{4a}
\left | \frac d {dh} s^k_l(h,g)_{\big\vert_{h=0}} \right |\leq C e^{\kappa |g|}
\end{equation}
for some $C>0$ and $\kappa \geq 1$. Thus we obtain the assertion of the proposition for $|\beta|=1$ , and  $b_i^j(m,g) =\sum _{k=1}^d \Theta^k_i(m) \dot s^k_j(0,g)$. Assume that the assertion holds for $|\beta|\leq N$. Differentiating \eqref{2} with  respect to $m_i$ then yields
the assertion for $|\beta|=N+1$.  
\end{proof}
\begin{corollary}
Assume  $\tilde m \in V_\R- S_\R$, and let $U$ be a neighbourhood of $\tilde m$ as in the preceeding proposition. Then, for arbitrary multiindices $\alpha$, $N\in \N$, and $g \in G$, the relations
\begin{equation}
\label{5}
\gd ^\alpha_\xi \pi(g) \phi_\xi(m)=\frac 1 {(1+|\xi|^2)^N} \sum \limits_{|\gamma| \leq |2 N} d^\gamma_{\alpha,N}(m,g) \pi(g) d\pi(X^\gamma)\phi_\xi(m), \qquad m \in U, 
\end{equation}
hold, where the coefficients $d^\gamma_{\alpha,N}(m,g) \in \Cinft(U \times G)$ are rational functions in the coordinates of $m$  satisfying the estimates
\begin{equation}
\label{6}
|\gd _m^\beta d^\gamma_{\alpha,N}(m,g)|\leq C_U e^{\kappa |g|}
\end{equation}
for some $C_U, \kappa>0$.
\end{corollary}
\begin{proof}
The key step in proving the corollary will be to express $(1+\xi^2)^N$ as a linear combination of derivatives of $\pi(g)\phi_\xi(m)$ with respect to $m$.
Since $\gd^\beta_m\phi_\xi(m)=i^{|\beta|} \xi^\beta \phi_\xi(m)$,  one computes
\begin{equation}
\label{6a}
(\gd_{m_j} \pi(g)\phi_\xi)(m)=\sum_{k=1}^{n} \frac { \gd(g^{-1}m)_k } {\gd m_j}  (\gd _{m_k} \phi_\xi)(g^{-1}m)=i\phi_\xi(g^{-1}m) \sum_{k=1}^n (g^{-1})_{kj} \xi_k,
\end{equation}
and  repeated differentiation leads to
\begin{align}
\label{6b}
\begin{split}
(\gd^\beta_m \pi(g) \phi_\xi)(m)&=i^{|\beta|} \phi_\xi(g^{-1}m)\Big ((^Tg^{-1}) \xi\Big )^\beta\\&=\phi_{(1-|\beta|)\xi}(g^{-1}m)\prod_\beta \big (\gd_{m_j} \pi(g)\phi_\xi(m)\big )^{\beta_j}.
\end{split}
\end{align}
Therefore, by first taking into account that  \eqref{6a} implies
\begin{displaymath}
\xi_j=-i\phi_{-\xi}(g^{-1}m)\sum_{k=1}^n g_{kj}(\gd_{m_k} \pi(g)\phi_\xi)(m),
\end{displaymath}
and then applying  \eqref{6b}, we obtain for  $(1+\xi^2)^N$ the expression 
\begin{displaymath}
(1+|\xi|^2)^N=\sum \limits ^N_{k=0} \left ( \begin{array}{c} N \\ k \end{array}
\right ) (\xi_1^2+\dots +\xi^2_n)^k =\phi_{-\xi}(g^{-1}m) \sum \limits_{|\beta| \leq 2N} c_\beta(g)  (\gd^\beta_m \pi(g)\phi_\xi)(m),
\end{displaymath}
where the  $c_\beta(g)$ are rational functions  in the matrix coefficients of $g$. The corollary now follows, since by the previous proposition, one  computes
\begin{align*}
(\gd ^\alpha_\xi \pi(g) \phi_\xi)(m)&=(g^{-1}m)^\alpha i^{|\alpha|}  \phi_\xi(g^{-1}m)\\ &=\frac 1 {(1+|\xi|^2)^N}(g^{-1}m)^\alpha i^{|\alpha|}\sum \limits_{|\beta| \leq 2N} c_\beta(g)(\gd ^\beta _m\pi(g) \phi_\xi)(m)\\ &=\frac 1 {(1+|\xi|^2)^N}(g^{-1}m)^\alpha i^{|\alpha|}\sum \limits_{|\beta| \leq 2N} c_\beta (g)\sum \limits_{|\gamma|\leq |\beta|} b^\gamma_\beta(m,g) (\pi(g)d\pi(X^\gamma)\phi_\xi)(m),
\end{align*}
where the functions $b^\gamma_\beta(m,g)$ are  rational expressions in the cordinate functions of $m$  satisfying the bounds \eqref{3}.
\end{proof}
We are now in position to prove the structure theorem.
\begin{proof}[Proof of Theorem \ref{thm:1}]
Fix  $\tilde
m\in V_\R- S_\R$, and assume that $U\subset V_\R- S_\R$ is a neighbourhood of  $\tilde m$ as in Proposition \ref{prop:1}. Then, by Proposition \ref{prop:1} and its corollary, one has 
\begin{align*}
\gd^\alpha_\xi \gd ^\beta_m \hat f(m,\xi)&=\int _G f(g) \gd
^\alpha_\xi \gd^\beta_m \pi(g) \phi_\xi(m) d_G(g)\\&=\frac 1 {(1+|\xi|^2)^N} \int_G f(g) 
\gd^\beta_m \Big [\sum \limits_{|\gamma|\leq 2N}
d^\gamma_{\alpha,N}(m,g)\pi(g) d\pi(X^\gamma)\phi _\xi (m)\Big ]
d_G(g)\\
&=\frac 1 {(1+|\xi|^2)^N} \int_G f(g)\sum \limits_{|\gamma|\leq 2N}\sum
_{\delta_1+\delta_2=\beta} \frac{\beta!}{\delta_1!\delta_2!}\gd
_m^{\delta_1}d^\gamma_{\alpha,N}(m,g)\\
&\cdot \sum_{|\epsilon|\leq|\delta_2|} b^\epsilon_{\delta_2}(m,g)
\pi(g)d\pi(X^\epsilon X^\gamma)\phi_\xi(m) d_G(g), \qquad m \in U.
\end{align*}
Next, for arbitrary $\phi\in \Cinft(V_\R)$ and $X\in \U$, 
\begin{equation}
\label{10}
\pi(g)d\pi(X)\phi(m)=dR(X) \phi_m(g), \qquad g \in G, \, m\in V_\R,
\end{equation}
where we set $\phi_m(g)=\pi(g)\phi(m)$. Indeed, 
\begin{align*}
\pi(g) d\pi(X_i) \phi(m)&=\lim _{h \to 0} h^{-1} [\phi(\e{-hX_i}
g^{-1}m) - \phi(g^{-1}m)]=\lim _{h \to 0} h^{-1} [\phi_m(g\e{hX_i}) - \phi_m(g)],
\end{align*}
so \eqref{10} is correct for $X\in \g$. Assuming that the assertion holds for $X^\gamma\in \U$ with $|\gamma|=N$, we obtain for arbitrary $X_i$
 \begin{align*}
\pi(g) d\pi(X_iX^\gamma) \phi(m)&=\lim _{h \to 0} h^{-1} [d\pi(X^\gamma)\phi(\e{-hX_i}
g^{-1}m) - d\pi(X^\gamma)\phi(g^{-1}m)]\\&=\lim _{h \to 0} h^{-1}
 [dR(X^\gamma)\phi_m(g\e{hX_i}) - dR(X^\gamma) \phi_m(g)]=dR(X_iX^\gamma)\phi_m(g),
\end{align*}
so that, by induction, we get \eqref{10} for arbitrary $X\in \U$. Now, integrating by parts according to \eqref{B} yields
 \begin{align*}
\gd^\alpha_\xi \gd ^\beta_m \hat f(m,\xi)& =\frac 1 {(1+|\xi|^2)^N}\sum \limits_{|\gamma|\leq 2N}\sum
_{\delta_1+\delta_2=\beta}  \sum_{|\epsilon|\leq|\delta_2|} \int_G
F^{\alpha,\beta,N}_{\delta_1,\delta_2,\gamma,\epsilon}(m,g) dR(X^\epsilon X^\gamma)\phi_{\xi,m}(g) d_G(g)\\&=\frac 1 {(1+|\xi|^2)^N}\sum \limits_{|\gamma|\leq 2N}\sum
_{\delta_1+\delta_2=\beta}  \sum_{|\epsilon|\leq|\delta_2|}
\sum_{\sigma_1+\sigma_2=\rho(\epsilon,\gamma)} \iota^{\sigma_1}
(-1)^{|\sigma_2|}\\&\cdot \int_G
dR(X^{\tilde \sigma_2})F^{\alpha,\beta,N}_{\delta_1,\delta_2,\gamma,\epsilon}(m,g) \pi(g)\phi_{\xi}(m) d_G(g),
\end{align*}
where we set
\begin{displaymath}
F^{\alpha,\beta,N}_{\delta_1,\delta_2,\gamma,\epsilon}(m,g)=\frac{\beta!}{\delta_1!\delta_2!}f(g)(\gd
_m^{\delta_1}d^\gamma_{\alpha,N})(m,g) b^\epsilon_{\delta_2}(m,g),
\end{displaymath}
and $X^{\rho(\epsilon,\gamma)}=X^\epsilon X^\gamma$. Note that,  for fixed  $m\in U$,
$F^{\alpha,\beta,N}_{\delta_1,\delta_2,\gamma,\epsilon}(m,g)\in \S(G)$, as a consequence of the estimates \eqref{3}
and \eqref{6}, so that integration by parts is legitimate. 
Let $\omega$ be a compact set in $V_\R$. For each point $\tilde
m\in \omega$, let  $U_{\tilde m}\subset V_\R- S_\R$ be a neighbourhood of $\tilde m$ as in proposition \ref{prop:1}. By Heine-Borel, $\omega$ can be  covered by finitely many neighbourhoods $U_{\tilde m}$, so that $\bigcup _{\tilde m \in \omega} U_{\tilde m}$
has  a finite subcovering. Fix $l\in \R$, and let $k$  be an integer $\leq l$. Assume that $\alpha,\, \beta$ are arbitrary multiindices. Setting   $N=|\alpha|+|k|$, and taking into account the above expression for $\gd^\alpha_\xi \gd ^\beta_m \hat f(m,\xi)$, one computes
\begin{equation*}
|\gd^\alpha_\xi \gd ^\beta_m \hat f(m,\xi)|\leq C_{\alpha,\beta,
 N,\omega} \eklm{\xi}^{-2N} \leq C_{\alpha,\beta,k,\omega}
 \eklm{\xi}^{k-|\alpha|}\leq C_{\alpha,\beta,l,\omega}
 \eklm{\xi}^{l-|\alpha|}, \qquad m \in \omega,
\end{equation*}
where $\eklm{\xi}=(1+|\xi|^2)^{1/2}$. Thus, $\hat
f(m,\xi)\in \Sym^l_{1,0}((V_\R- S_\R) \times V_\R^\ast)$ for all $l \in \R$, yielding $a_f(m,\xi)\in \Syms((V_\R- S_\R) \times V_\R^\ast)$.
Now, for $\phi\in \CT(V_\R- S_\R)$, Fourier transformation gives
\begin{align*}
\pi(f)\phi(m)&= \int _G f(g)\pi(g)\phi (m) \d_G(g)=\int_G f(g)
\left (\int e^{i(g^{-1}m)\cdot \xi} \hat \phi(\xi) \dbar \xi\right )
\d_G(g)\\&=\int \hat f(m,\xi) \hat \phi(\xi) \dbar(\xi) =\int \int 
e^{i(m-m')\cdot \xi} a_f(m,\xi) \phi(m') \d m' \, \dbar \xi,
\end{align*}
where the occuring integrals are absolutely convergent, so that we can interchange the order of integration. This proves Theorem \ref{thm:1}. 
\end{proof}

\section{Fixed point singular sets}
\label{sec:V}

Let $(G_\C,\rho,V_\C)$ be a prehomogeneous vector space with singular set $S_\C$. 
In this section, we restrict our attention to the case where $S_\R=S_\C\cap V_\R$ coincides with the set of fixed points of the $G_\R$-action. If $G_\R$ is reductive, zero constitutes the only closed orbit contained in $S_\R$, and hence the only fixed point. This will be the case we will be mainly concerned with.
If $m\in V_\R$ is a fixed point, the symbol of $\pi(f)$ at $m$ is simply given by $a_f(m,\xi)=\int _G f(g) \d_G(g)$. Thus, if $S_\R$ coincides with the set of fixed points $F_\R$ of the underlying $G_\R$-action, one has 
\begin{equation}
\label{11}
\begin{split}
a_f(m,\xi)\in \left \{\begin{array}{ll}
 \Syms((V_\R)_{m'}\times V_\R^\ast) & \text{if } m\in
V_\R- F_\R,\\ \Sym^0_{1,0}((V_\R)_{m'}\times V_\R^\ast) &
\text{if } m \in F_\R.\qquad\end{array}\right.
\end{split}
\end{equation}
In this situation, Theorem \ref{thm:1} can be restated as follows.
\begin{theorem}
\label{thm:2}
Let $(G_\C,\rho,V_C)$ be a prehomogeneous vector space defined over $\R$, and assume that $S_\R$ is equal to the set $F_\R$ of fixed points of the underlying $G_\R$-action. Then $\pi(f):\CT(V_\R)\rightarrow\D'(V_\R)$ is given by the family of oscillatory integrals
\begin{equation*}
(\pi(f)\phi)(m)=\int e^{i(m-{m'} )\cdot \xi} a_f(m,\xi) \phi({m'}) d{m'} \dbar \xi,
\end{equation*}
where 
\begin{equation*}
a_f(m,\xi)= e^{-im\cdot \xi} \hat f(m,\xi) \in \Sym^0_{1,0} ( (V_\R)_{m'} \times V_\R^\ast), \qquad m \in V_\R;
\end{equation*}
for  $m\notin F_\R$, one has $a_f(m,\xi)\in \Syms ((V_\R)_{m'}\times
V_\R^\ast)$. Its Schwartz kernel  $\K_f\in \D'(V_\R\times V_\R)$ is given by the continuous family of oscillatory integrals
\begin{equation*}
V_\R \ni m \mapsto \K_{f,m} = \int e^{i( m-\cdot)\cdot  \xi} a_f(m,\xi) \dbar \xi.
\end{equation*}
Furthermore, one has $a_f(m,\xi) \in \Syms((V_\R)_m-F_\R)\times (V_\R)_{m'}\times V_\R^\ast)$, so that $\pi(f)$, as a mapping from $\CT(V_\R)$ to   $\D'(V_\R-F_\R)$, represents a Fourier integral operator. In particular, $\pi(f)\in \L^{-\infty} (V_\R-F_\R)$.
\end{theorem}
\begin{example}
\label{ex:1}
Consider the simplest prehomogeneous vector space, $(\C^\ast,\C)$, where $\C^\ast$ acts by multiplication on $\C$. Then $S_\C=S_\R=\mklm{0}$, and $p(m)=m$ is a relative invariant corresponding to $\chi(g)=g$, $g \in \C^\ast$. One has $V_\R-S_R=\R^\ast=V_1\cup V_2$ with $V_1=R^\ast_+$, $V_2=\R^\ast_-$. The fundamental theorem yields in this case the relations
\begin{equation}
\label{12}
  \int_{-\infty}^\infty |\xi|^{s-1} \hat \phi(\xi) \d \xi =(2\pi)^{-s} \Gamma(s) \, 2 \cos \frac{\pi s} 2 \int_{-\infty}^\infty |m|^{-s} \phi(m) \d m, \qquad \phi \in \S(\R),
\end{equation}
see \cite{kimura}, Proposition 4.21.
\end{example}
\begin{example}
  Let $B$ be a positive definite real symmetric matrix of degree $n\geq 3$, and consider the positive definite quadratic forms $p(m)=m^t B m$,  $p^\ast(m)=m^tB^{-1} m$ on $\C^n$. We define $\SO(n,B)=\mklm{X \in \SL(n,\C): X^t BX=B}$. Setting $\rho(\alpha,g) m =\alpha g m$, $\rho^\ast(\alpha,g)=\alpha^{-1} (g^{-1})^t m$, where $\alpha \in \C^\ast$, $g \in \SO(n,B)$, one obtains actions of $G_\C=\GL(1,\C) \times \SO(n,B)$ on $V_\C=V_\C^\ast=\C^n$. Again,  $S_\R=\mklm{0}$, and by the fundamental theorem one obtains
  \begin{equation}
    \label{12a}
\int_{\R^n-\mklm{0}} |p^\ast(\xi)|^{s-\frac n2} \hat \phi(\xi) \d \xi= \pi^{\frac n2 -s} \sqrt{\det B} \frac {\Gamma(s)}{\Gamma(\frac n2 -s)} \int _{\R^n-\mklm{0}} |p(m)|^{-s} \phi(m) \d m,
  \end{equation}
where $\phi \in \S(\R^n)$, see \cite{kimura}, Proposition 4.22.
\end{example}

Returning to the situation of the previous theorem, note that on  $({V_\R}- F_\R)\times {V_\R}$, the Schwartz kernel  $\K_f$ of $\pi(f)$ is given by the function
\begin{equation*}
K_f(m,{m'})=\int e^{i(m-{m'})\cdot \xi} a_f(m,\xi) \dbar \xi \in \Cinft(({V_\R}- F_\R)\times {V_\R}),
\end{equation*}
while, as distributions,
$\K_{f,m} \to \K_{f,x_F}=\int_G f(g) d_G(g) \cdot \delta_{x_F}$ as $m\to x_F \in F_\R$. In this way, the orbit structure of the underlying $G_\R$-action is reflected in the singular behaviour of the Schwartz kernel of 
$\pi(f):\CT({V_\R}) \rightarrow \D'({V_\R})$ at $F_\R$. In what follows, we will assume that $S_\R=\mklm{0}$. In order to get a better understanding of  the Schwartz kernel  $\K_f$
of   $\pi(f)$, and, in particular, of its restriction to the diagonal, we define the auxiliary symbol
\begin{equation*}
\tilde a_f(m,\xi)=e^{-i\frac {m \cdot \xi}{|m|}}\hat f(m/|m|,\xi)=a_f(m/|m|,\xi), \qquad m \not=0.
\end{equation*}
By introducing the inverse Fourier transform
\begin{displaymath}
\tilde A_f(m,m')= \int e^{i {m'} \cdot\xi } \tilde a_f(m,\xi) \dbar \xi
\end{displaymath}
 of $\tilde a(m,\xi)$, and since $a_f(\alpha m,\xi)=a_f(m,\alpha \xi)$, $\alpha \in\R$, 
we obtain the relation
\begin{align*}
K_f(m,{m'})&=\int e^{i(m-{m'})\cdot\xi} \tilde a_f(m,|m|\xi) \dbar
\xi= \frac 1 {|m|^n} \tilde A_f\Big (m, \frac{m-{m'}}{|m|}\Big), \quad m \not=0.
\end{align*}
\begin{lemma}
\label{lem:2}
$K_f(m,{m'})\in \L^1_{\rm{loc}}(\R^{2n})\cap \Cinft(\R^{2n}-\mklm{0})$.
\end{lemma}
\begin{proof}
By the structure theorem, $\tilde a_f(m,\xi)\in  \Syms({V_\R}-\mklm{0},{V_\R})$, which implies that $\tilde A_f(m,{m'})$ vanishes to all orders as $|{m'}|\to \infty$. Since $\tilde A_f(m,{m'})$ only depends on the direction of $m$, $\tilde A_f(m,(m-{m'})/|m|)$ goes to zero to infinite order as $|m|\to 0$, provided  ${m'}\not=0$. By setting $K_f(m,{m'})=0$ for   $m=0,\, {m'}\not=0$, we therefore obtain $K_f(m,{m'})\in\Cinft(\R^{2n}- \mklm{0})$.
Let $L$ be a compact set in $\R^{2n}$. Then, for arbitrary $N$,  
\begin{displaymath}
\modulus{\int _L K_f(m,{m'}) \d m \d {m'} }\leq C_N \int _{L-
\mklm{0}} \frac {|m|^{2N-n}}{(|m|^2 +|m-{m'}|^2)^N} \d m \d {m'}. 
\end{displaymath}
For $N=n/2$, the last integral reads 
\begin{displaymath}
\int _{L- \mklm{0}} \frac { \d m \d
{m'}}{(|m|^2 +|m-{m'}|^2)^{n/2}}=\int\limits_0^R \int \limits_{S^{2n-1}} \frac {r^{2n-1} \d r \d \omega}{r^n\sum_{i=1}^n \omega_i^2+(\omega_i-\omega_{n+i})^2}, 
\end{displaymath}
where we introduced in $V_\R\times V_\R\simeq\R^{2n}$ the polar coordinates $r\in \R^+$ and $\omega \in S^{2n-1} \subset \R^{2n}$, and assumed that $L$ is contained in a sphere of radius $R$. Since   $\sum_{i=1}^n \omega_i^2+(\omega_i-\omega_{n+i})^2>0$ is a real valued, continuous function on the $(2n-1)$-dimensional sphere which is  bounded from below, the infimum is adopted. Hence, there is a strictly positive number $\kappa$ such that $\sum_{i=1}^n \omega_i^2+(\omega_i-\omega_{n+i})^2\geq \kappa$. This proves the lemma.
\end{proof}
The last lemma implies that $\K_f$ is given by the locally integrable function $K_f(m,m')$. Let us examine its restriction to the diagonal.  By the structure theorem,
\begin{equation*}
  k_f(m)=K_f(m,m)=\int a_f(m,\xi) \dbar \, \xi=\frac 1 {|m|^n} \tilde A_f (m,0) \in \Cinft(V_\R-\mklm{0}).
\end{equation*}
As a consequence of Theorem \ref{thm:A}, we then have the following proposition.
\begin{proposition}
Let $(G_\C,\rho,V_\C)$ be a reductive prehomogeneous vector space, and $S_\R=\mklm{0}$. Set $k_{f,s}(m)=|m|^s \tilde A_f(m,0)$. Then, for  $\phi \in \CT(\R^n)$,  the integrals
\begin{equation*}
 \int _{\R^n-\mklm{0}} k_{f,s}(m) \phi(m) \d m , 
\end{equation*}
can be continued to meromorphic functions on the whole $s$-plane. For $\phi\in\CT(\R^n-\mklm{0})$, they satisfy the functional equations
\begin{equation*}
\int_{\R^n-\mklm{0}} |\xi|^{2(s-\frac n2)} \widehat{\tilde A_f(\cdot,0)\phi}(\xi) \d \xi =\pi^{\frac n2 -2s} \frac {\Gamma(s)}{ \Gamma\Big ( \frac n2 -s\Big )}\int _{\R^n-\mklm{0}} k_{f,-2s}(m) \phi(m) \d m.
\end{equation*}
\end{proposition}
\begin{proof}
Put $p(m)=m_1^2+\dots+m_n^2=|m|^2$, and assume $\Re s > -n/2$. One then computes
\begin{align*}
  \int_{\R^n-\mklm{0}}k_{f,2s}(m) \phi(m) \d m & =\int_{\R^n-\mklm{0}} |p(m)|^s \tilde A_f(m,0) \phi(m) \d m  \\ &=\int_{S^{n-1}} \tilde  A_f(\omega,0) \int _0^\infty \phi(r\omega) r^{2s+n-1} dr d\omega,
\end{align*}
everything in sight being absolutely convergent. Consider now the $n$-dimensional Mellin transform
\begin{equation*}
  M_s: \S(\R^n) \longrightarrow \Cinft(S^{n-1}), \qquad M_s (\phi)(\omega) = \int _0^\infty \phi(r\omega) r^s dr, \qquad \Re s > -1,
\end{equation*}
which can be continued meromorphically in $s$ with simple poles at $s=-1,-2,\dots$. Since, by the above computation, 
\begin{equation*}
   \int_{\R^n-\mklm{0}}k_{f,2s}(m) \phi(m) \d m=\int_{S^{n-1}} \tilde A_f(\omega,0) M_{2s+n-1} (\phi) (\omega) d\omega,
\end{equation*}
we obtain the first assertion. The second one follows directly from the relations \eqref{12a} with respect to the quadratic form $p(m)$. 
\end{proof}
Now, although the integrals $\int _{\R^n-\mklm{0}} k_{f}(m) \phi(m) \d m $ become singular for general $\phi \in \CT(\R^n)$, it is nevertheless possible to extend $k_f$, as a distribution, to $\R^n$. First note that, since  $K_f(tm,tm)=t^{-n}K_f(m,m)$ on $V_\R-\mklm{0}$,  where $t>0$,  $k_f(m)$ defines a homogeneous distribution on $V_\R-\mklm{0}$ of degree $-n$ according to
\begin{equation*}
  \eklm{k_f,\phi}=\int _{|\omega|=1} \int^\infty_0 k_f(w) r^{-1} \phi(r\omega) \,dr \,d\omega.
\end{equation*}
By \cite{hoermanderI}, Theorem 3.2.4, an extension of $k_f$ to $V_{\R}$ can then be constructed as follows. Let $s$ be a complex number with $\Re s >-1$, and consider on $\R$ the locally integrable function
\begin{equation*}
  t_+^s=t^s, \quad \text{if  }t >0, \qquad t_+^s=0 \quad \text {if } t \leq 0.
\end{equation*}
It is homogeneous of degree $s$, and, as a distribution, can be extended to arbitrary complex values of $s$ by analytic continuation, except when $s=-1,-2,\dots$. In that case one defines for $\phi \in \CT(\R)$
\begin{equation*}
  t^{-k}_+(\phi)=-\frac 1 {(k-1)!} \int^\infty_0 (\log t) \phi^{(k)}(t) dt + \frac 1 {(k-1)!}\phi^{(k-1)} (0) \Big(\sum_1^{k-1}\Big )j^{-1}.
\end{equation*}
For $s\not=-1,-2,\dots$, $t^s_+$ is then a homogeneous distribution on $\R$.
Now, for $m \not=0$, and $\phi\in \CT(V_\R)$, set $R_s\phi(m)=t^{s+n-1}_+(\phi(tm))$, and let $\psi \in \CT(V_\R-\mklm{0})$ satisfy
\begin{equation*}
  \int _0^\infty \psi(tm) \frac {dt} t=1, \quad m \not=0.
\end{equation*}
Then, the extension of $k_f$ is defined as 
\begin{equation}
\label{13}
  \dot {k_f}(\phi)=\eklm{k_f, \psi R_{-n}\phi}, \qquad \phi \in \CT(V_\R).
\end{equation}
Summing up, we have shown the following proposition.
\begin{proposition}
\label{prop:2}
Let $(G_\C,\rho,V_\C)$ be a prehomogeneous vector space, and $S_\R=\mklm{0}$. Then, the restriction of the kernel of $\pi(f)\in \L^{-\infty}(V_\R-\mklm{0})$ to the diagonal is given by the homogeneous distribution $k_f=K_f(m,m)\in \Cinft(V_\R-\mklm{0})$ of degree $-n$, which has an extension to $V_\R$ given by \eqref{13}.
\end{proposition}
\begin{remark} The extension $\dot {k_f}$ is unique up to a linear combination of $\delta$-distributions at zero, and is no longer homogeneous.
 One could also have defined $\dot {k_f}$ as the integral over the unit sphere of $k_f R_{-n}\phi$ according to
  \begin{equation*}
    \dot {k_f}(\phi)=\int_{|\omega|=1} k_f(\omega)R_{-n}\phi(\omega) \d\omega, \qquad \phi \in \CT(V_\R).
  \end{equation*} 
Each of the extensions $\dot k_f$ could then be regarded as a trace of $\pi(f)$.
\end{remark}

\section{Totally characteristic pseudodifferential operators on prehomogeneous vector spaces}

Let $(G_\C,\rho,V_\C)$ be a reductive prehomogeneous vector space defined over $\R$, and assume that the singular set $S_\C$ is an irreducible hypersurface.  Denote by $p$ the corresponding relative invariant with character $\chi$ such that \eqref{E} is satisfied, and  by $(\pi,\Cvan(V_\R))$ the left regular representation of $G_\R$ on the Banach space $\Cvan(V_\R)$. As already explained in Section \ref{sec:IV}, for $f\in \S(G)$, the restriction of the operator  $\pi(f)$ to $\CT(V_\R)$ defines a continuous linear operator 
\begin{equation*}
  \pi(f):\CT(V_\R) \rightarrow \Cinft(V_\R).
\end{equation*}
Now, if $V_\R-S_\R=V_1\cup \dots \cup V_l$ denotes the decomposition of the open $G_\R$-orbit into its connected components, on each $V_i$, $\pi(f)\phi$ only depends on the restriction of $\phi\in\Cvan(V_\R)$ to $V_i$, so that one naturally obtains the continuous linear operators 
\begin{equation*}
  \pi(f)_{|V_i}: \CT(V_i)\longrightarrow  \Cinft(V_i), 
\end{equation*}
which are elements in $\L^{-\infty}(V_i)$, by the structure theorem. Throughout this section, we will make the following assumption.
\begin{assumption}
$S_\R=S_\C \cap V_\R$ is an irreducible hypersurface.  
\end{assumption}
Denote the set of non-regular points of $S_\R$ by $S_\R^{\rm{sing}}$. We will then show that the restrictions of $\pi(f)$ to ${\overline V_i -S_\R^{\rm{sing}}}$ are totally characteristic pseudodifferential operators in the class $\L^{-\infty}_b$. Let us introduce local coordinates in $V_\R-S_\R^{\rm{sing}}$. Since the fibers of the categorical quotient $p:V_\R \rightarrow V_\R/G_\R$ are, apart from the exceptional divisor $S_\R$, smooth affine varieties, we have $\grad p(m) \not=0$ on $V_\R-S_\R^{\rm{sing}}$. Defining  the open subsets
\begin{equation*}
  (V_\R-S_\R^{\rm{sing}})_j=\mklm {m \in V_\R-S_\R^{\rm{sing}}: (\gd/\gd_{m_j} p)(m)\not=0},
\end{equation*}
and the coordinates 
\begin{equation*}
  \kappa(m)=(x_1(m),\dots,x_n(m)), \qquad x_1(m)=p(m), \quad x_2(m)=m_{l_1},\dots, x_n(m)=m_{l_{n-1}}, 
\end{equation*}
on $(V_\R-S_\R^{\rm{sing}})_j$, where $\mklm{m_{l_1},\dots, m_{l_{n-1}}}\cup \mklm{j}$, one gets
\begin{equation}
  \det \left (\frac {\gd x_i}{\gd m_k}(m) \right )_{ik} \not= 0  \qquad \text{for } m \in (V_\R-S_\R^{\rm{sing}})_j.
\end{equation}
By the inverse function theorem, there exists for every $m \in  (V_\R-S_\R^{\rm{sing}})_j$ an open neighbourhood $W_j\subset (V_\R-S_\R^{\rm{sing}})_j$ of $m$ such that $\kappa: W_j \rightarrow \tilde W_j=\kappa (W_j) \subset \R^n$
becomes a diffeomorphism, and we obtain a covering of $V_\R-S_\R^{\rm{sing}}$ by charts.  We define now the following closed subgroups of $G_\C$. Let
\begin{equation*}
  G_\C(p)=\mklm {g \in G_\C: p(gm )=p(m)}=\mklm{g \in G_\C: \chi(g) =1}, \qquad G_\R(p)= G_\C(p) \cap G_\R.
\end{equation*}
Clearly, $G_\R(p)$ acts transitively on each generic fiber of the categorical quotient  $p:V_\R \rightarrow V_\R/G_\R$, and its Lie algebra is given by $\g(p)=\mklm {X\in \g: d\chi(X) =0}$, where $d\chi$ denotes the infinitesimal representation corresponding to $\chi$. There always exists a  $Y\in \g$ such that $d\chi(Y)\not=0$; otherwise, $V_\R-S_\R$ would decompose into an infinite number of $G$-orbits. For the following considerations, we will make the basic assumption that $G_\R(p)$ acts transitively on $S_\R^{\rm{reg}}$. Let us formulate this condition as follows.
\begin{assumption}
For all $m \in  V_\R-S_\R^{\rm{sing}}$, there exists an open neighbourhood $Z \subset V_\R-S_\R^{\rm{sing}}$  of $m$ and elements $X_1,\dots, X_{n-1} \in \g(p)$ such that $m'\to (\tilde X_1(m'),\dots, \tilde X_{n-1}(m'))$ defines a section in $T (p^{-1}(c) \cap Z)$ for each $c \in \R$. 
\end{assumption}
Here we wrote $\tilde X(m)=\frac d {dh}( \rho(\e{hX})m)_{|h=0}$ for the fundamental vector field of the underlying $G_\R$-action on $V_\R$ induced by $X\in \g$. Fix $m \in  (V_\R-S_\R^{\rm{sing}})_j$, and let $W_j$ and $Z$ be neighbourhoods of $m$ as specified above. Let $U,U_1$ be open sets containing $m$ such that $ U\subset {U}_1 \subset W_j\cap Z$, and assume that $U$ is chosen in such a way that the set $\mklm{g \in G_\R: g {U} \subset  {U}_1}$ acts transitively on the $G_\R$-orbits of $U$. The so defined  sets $U$ constitute an open covering of $V_\R-S_\R^{\rm{sing}}$. By choosing a locally finite subcovering, we therefore obtain an atlas $(\kappa_\gamma,U^\gamma)$  of $V_\R-S_\R^{\rm{sing}}$ with the following properties:

\bigskip

\noindent
(i) for $m \in U^\gamma$, one has $\kappa_\gamma(m)=(x_1,\dots,x_n)=(p(m),m_{l_1},\dots, m_{l_{n-1}})$, the $m_{l_j}$ being given by the prescription $\mklm{m_{l_1},\dots,m_{l_{n-1}},j}=\mklm{1,\dots,n}$ for some $j=j(\gamma)$;

\noindent
(ii) for each $U^\gamma$, there exist open sets $G^\gamma\subset G^\gamma_1\subset G_\R$, stable under inverse, acting transitively on the $G_\R$-orbits of $U^\gamma$;

\noindent
(iii) for $m \in \bigcup _{g \in G^\gamma_1} g U^\gamma$,  one has $\frac \gd {\gd m_j} p(m)\not=0$;

\noindent
(iv)   there exist $X_1,\dots, X_{n-1}\in \g(p)$, such that  $m'\to (\tilde X_1(m'),\dots, \tilde X_{n-1}(m'))$ defines a section in $T (p^{-1}(c) \cap \bigcup _{g \in G^\gamma_1} g U^\gamma)$ for each $c \in \R$.

\begin{example}
\label{ex:3}
  As the generic case in the theory of prehomogeneous vector spaces, consider an indefinite quadratic form $p$ with signature $(q,n-q)$, $1 \leq q \leq n-1$, on $\R$. Then, by Sylvester's law, $p(m)= B^t \1_{q,n-q} m^t B m$, for some $B \in \GL(n,\R)$. Define  the orthogonal, resp.\ special orthogonal, group of $p$ by $\O(p)=\mklm{g \in \GL(n,\C):p(gm)=p(m)}$, resp.\ $\SO(p)=\mklm{g \in \SL(n,\C):p(gm)=p(m)}$,  and introduce on $V_\C=\C^n$ an action  of $G_\C=\GL(1,\C) \times \SO(p)$ by setting $\rho(\alpha,g) m=\alpha g m$. Then $p$ is an irreducible relative invariant of $(G_\C,\rho,V_\C)$, and $V_\R-S_\R=V_+\cup V_-$, where $V_\pm=\mklm{ m \in V_R:\pm p(m) >0}$. By the theorem of Witt, $p^{-1}(0)-\mklm{0}$ is $\O(p)$-homogeneous, and $\SO(p)$-homogeneous for $n\geq 3$. In this case, the fundamental theorem yields the functional equations
  \begin{equation}
\label{22}
    \left ( \begin{array}{c} \int_{V_+^\ast} |p^\ast(\xi)|^{s-\frac n2} \hat \phi (\xi) \d \xi \\\int_{V_-^\ast} |p^\ast(\xi)|^{s-\frac n2} \hat \phi (\xi) \d \xi \end{array}\right ) = C(s)  \left ( \begin{array}{c} \int_{V_+} |p(m)|^{-s}      \phi (m) \d m \\\int_{V_-} |p(m)|^{-s} \phi (m) \d m \end{array}\right ),
  \end{equation}
where $\phi \in \S(V_\R)$, and 
\begin{equation*}
  C(s)= \Gamma\Big (s+1-\frac n2\Big ) \Gamma(s) {|\det B|} \pi^{-2s+\frac n2 -1} \left (
  \begin{array}{cc} -\sin \pi( s -\frac q2) & \sin \frac {\pi q}2 \\ \sin \frac {\pi(n-q)}2 & -\sin \pi( s -\frac {n-q}2 )\end{array}\right ),
\end{equation*}
see \cite{kimura}, Proposition 4.27.
\end{example}
\begin{example}
As an example of a non-regular prehomogeneous vector space, consider the subgroup $G_\C$ of upper triangular $2\times2$-matrices in $\GL(2,\C)$ acting regularly on  $V_\C=\C^2$. In this case, $S_\C=\mklm{m \in \C^2: m_2=0}$.
\end{example}

According to the structure theorem of the foregoing section,  $\pi(f)$ is a pseudodifferential operator on  $V_\R-S_\R$ with smooth kernel. 
In the present section, we will be concerned with the action of  $\pi(f)$ as an operator on $V_\R-S_\R^{\rm{sing}}$. With respect to the atlas $(\kappa_\gamma, U^\gamma)$, we define for $u \in \CT(U^{\gamma})$ the operators $A_f^{\gamma} u = [\pi(f)_{|U^\gamma} ( u\circ \kappa_{\gamma})] \circ \kappa_{\gamma}^{-1}$, see Section \ref{sec:III}. Explicitely, one has
  \begin{equation*}
    A_f^{\gamma} u(x)=\int_{G} f(g) \pi(g) (u \circ \kappa_{\gamma})(\kappa^{-1}_\gamma(x)) d_G(g), \qquad x \in \tilde U^\gamma.
  \end{equation*}
Let $c_\gamma\in \Cinft(G_\R)$ be a smooth function on $G_\R$ with support in $G^\gamma_1$ satisfying $c_\gamma \equiv 1$ on $G^\gamma$. Then, $A_f^\gamma$ can be written as
 \begin{equation*}
    A_f^{\gamma} u(x)=\int_{G} f(g)  (u ( \kappa^{\gamma}_g)(x)) c_\gamma(g) d_G(g), \qquad x \in \tilde U^\gamma,
  \end{equation*}
where we put $\kappa^\gamma_g=\kappa_\gamma \circ g^{-1} \circ \kappa^{-1}_\gamma$. We define 
\begin{gather}
  \hat f _\gamma(x,\xi)= \int _{G_\R} f(g) e^{i\kappa_g^\gamma(x)\cdot \xi} c_\gamma(g) d_G(g), \qquad 
 a^\gamma_f(x,\xi)=e^{-ix\cdot\xi} \hat f_\gamma(x,\xi).
\end{gather}
By Lebesgue's bounded convergence theorem, we can differentiate under the integral sign, yielding $\hat f _\gamma(x,\xi)$, $ a^\gamma_f(x,\xi)\in \Cinft(\tilde U^\gamma\times \R^n_\xi)$. Set 
\begin{equation*}
  T_x=\left (\begin{array}{cc} x_1 & 0 \\ 0 & \1_{n-1} \end{array} \right ), \qquad x \in \tilde U^\gamma,
\end{equation*}
and define the symbol 
\begin{equation*}
 \tilde a^\gamma_f(x,\xi)=a_f^\gamma(x,T_x^{-1}\xi), \qquad x_1 \not= 0.
\end{equation*}
Note that, by equation \eqref{E},  $p(g\kappa^{-1}_\gamma(x))=\chi (g) p(\kappa^{-1}_\gamma(x)) =\chi(g) x_1$, so that
\begin{align*}
  T_x^{-1} \kappa ^\gamma_g(x)&=T_x^{-1} (x_1(g^{-1} \kappa ^{-1}_\gamma (x)), \dots, x_n(g^{-1} \kappa ^{-1}_\gamma (x)))=(\chi(g^{-1}), m_{l_1}(g^{-1} \kappa ^{-1}_\gamma (x)),\dots),
\end{align*}
implying that $\tilde a_f^\gamma(x,\xi)$ can be expressed explicitely by
\begin{equation}
\label{21}
 \tilde a_f^\gamma(x,\xi)=e^{-i(1,x_2,\dots,x_n)\cdot \xi} \int_{G_\R} f(g) e^{i(\chi(g^{-1}),m_{l_1}(g^{-1} \kappa^{-1}_\gamma(x)),\dots) \cdot \xi} c_\gamma(g) d_G(g),
\end{equation}
yielding  $\tilde a_f^\gamma(x,\xi)\in \Cinft(\tilde U^\gamma\times \R^n_\xi)$.  Put 
\begin{equation*}
 \tilde U^\gamma_+=\mklm{ x \in \tilde U^\gamma: x_1 \geq 0}, \qquad \tilde U^\gamma_-=\mklm{ x \in \tilde U^\gamma: x_1 \leq 0}.
\end{equation*}
 The main result of this section will consist in proving the following theorem. Here $(\kappa_\gamma,U^\gamma)$ denotes the atlas of $V_\R-S_\R^{\rm{sing}}$ constructed above.
\begin{theorem}
\label{thm:3}
Let $(G_\C,\rho,V_\R)$ be a reductive prehomogeneous vector space whose singular set $S_\C$ is an irreducible hypersurface, and let Assumptions $1$ and $2$ be satisfied. Then, for $f\in \S(G)$, $\pi(f)$ is locally given by the operators
\begin{gather}
\label{20}
  A^\gamma_fu(x)= \int e ^{i (x-y) \cdot \xi} a_f^\gamma(x,\xi)u(y) \d y \,\dbar\xi, \qquad u \in \CT(\tilde U^\gamma),
\end{gather}
where $a_f^\gamma(x,\xi)=\tilde a_f^\gamma(x,T_x\xi)$, and $\tilde a_f^\gamma(x,\xi) \in \Symsl(\tilde U^\gamma \times \R^n_\xi)$ is given by \eqref{21}. In particular, $A_f^\gamma \in \L^{-\infty}_b({\tilde U^\gamma_\pm})$.
\end{theorem}

Note that the symbols $a_f^\gamma(\kappa_\gamma(m),\xi)$ are defined even in the case when  the coordinates $\kappa_\gamma$ become singular. The expressions \eqref{20} are then still absolutely convergent, so that $\pi(f)$ is indeed given globally by the operators $A_f^\gamma$.
We will divide the proof of  Theorem \ref{thm:3} in several parts. To begin with, let us first state a technical lemma.
\begin{lemma}
  Let $\tilde \phi_{\xi,x}^\gamma (g)=\phi_{T^{-1}_x\xi}(\kappa^\gamma_g(x))$. Then, for $X \in \g$, $g \in G_\R$, and $\bar \chi(g)= \chi(g^{-1})$ one has 
  \begin{equation*}
    dR(X) \tilde \phi^\gamma_{\xi,x}(g) = i \tilde \phi _{\xi,x}^\gamma (g) [ \xi_1 d\chi(-X) \bar \chi(g) +\sum _{j=2}^n \xi_j m_{l_{j-1}} (-X g^{-1} \kappa ^{-1}_\gamma(x))].
  \end{equation*}
\end{lemma}
\begin{proof}
  First, one computes
  \begin{align*}
    dR(X) \tilde \phi^\gamma_{\xi,x}(g)&=\frac d {dh} {e ^{i (T_x^{-1}\kappa^\gamma_{g \e{hX}}(x)) \cdot \xi}}_{|_{h=0}}=i \tilde \phi_{\xi,x}^\gamma(g) \sum _{j=1}^n \xi_j \frac d {dh}{ \Big (T_x^{-1} \kappa^\gamma_{g\e{hX}} (x)\Big ) _j}_{|_{h=0}}\\&=i\tilde \phi^\gamma_{\xi,x}(g) \Big [ \xi_1 dR(X) \bar \chi(g) + \sum ^n_{j=2} \xi_j d R(X) m_{l_{j-1},\kappa^{-1}(x)}(g)\Big ].
  \end{align*}
By noting that
\begin{equation*}
  dR(X) \chi (e)=\frac d {dh} \chi(\e{hX})_{|_{h=0}}=d\chi(X),
\end{equation*}
we obtain $dR(X) \bar \chi (g) =d \chi(-X) \bar \chi(g)$. Similarly, one verifies the relation
$ dR(X) m_{l_j,\kappa^{-1}(x)}(g) =-m_{l_j}(Xg^{-1}\kappa^{-1}(x))$, thus obtaining the assertion.
\end{proof}
\begin{remark} For $X \in \g$, one has
  \begin{equation*}
    d\pi(-X) p(m)=\frac d {dh} p(\rho(\e{hX}) m)_{|_{h=0}}=d\chi(X) p(m),
  \end{equation*}
 as well as
 \begin{equation*}
   d\chi(X) p(m) = (X m)\cdot \grad p(m).
 \end{equation*}
Thus, $d\chi(X)$ is a measure for the transversality of the fundamental vector field $\tilde X(m)=X m$ with respect to the foliation
\begin{equation*}
  V_\R=\cup_{c \in \R} (p^{-1}(c) \cap V_\R).
\end{equation*}
\end{remark}
As a first step towards Theorem \ref{thm:3}, we prove the following 
\begin{proposition}
  \label{prop:3}
$\tilde a^\gamma_f(x,\xi) \in \Syms(\tilde U^\gamma \times \R^n_\xi)$.
\end{proposition}
\begin{proof}
  As already noted, $\tilde a ^\gamma_f(x,\xi) \in \Cinft(\tilde U^\gamma \times \R^n_\xi)$. While differentiation with respect to $\xi$ does not alter the growth properties of $\tilde a_f^\gamma(x,\xi)$, differentiation with respect to $x$ yields additional powers in $\xi$. Let $Y\in\g$ and $X_1,\dots,X_{n-1}\in \g(p)$  be given in such a way that $d\chi(Y) \not=0$, and condition (iv) is satisfied. By the previous lemma,
  \begin{equation*}
    \left ( 
    \begin{array}{c} dR(X_1) \tilde \phi^\gamma_{\xi,x}(g) \\ \vdots \\  dR(X_{n-1}) \tilde \phi^\gamma_{\xi,x}(g) \\ dR(Y) \tilde \phi_{\xi,x}^\gamma(g) \end{array} \right )=i \tilde \phi_{\xi,x}^\gamma(g) \left ( 
    \begin{array}{ccc} d\chi(-X_1) \bar \chi(g) & m_{l_1} (-X_1 g^{-1}\kappa^{-1}_\gamma(x)) &\dots \\
   \vdots & \vdots & \ddots \\ d\chi(-X_{n-1}) \bar \chi(g) & m_{l_1} ( -X_{n-1} g^{-1} \kappa^{-1} \gamma(x)) &\dots  \\ d\chi(Y) \bar \chi(g)& m_{l_1}(-Yg^{-1}\kappa^{-1}_\gamma (x)) & \dots      
    \end{array} \right ) \xi
  \end{equation*}
Writing $\Gamma(x,g)$ for the matrix appearing on the right hand side of the last equation, we obtain
\begin{equation*}
  \det \Gamma(x,g) = d\chi(-Y)\bar \chi (g) \left | 
  \begin{array} {ccc}
  m_{l_1} (-X_1 g^{-1}\kappa^{-1}_\gamma(x)) & \dots & m_{l_{n-1}} (-X_1 g^{-1}\kappa^{-1}_\gamma(x)) \\ \vdots &\ddots & \vdots \\ m_{l_1} (-X_{n-1} g^{-1}\kappa^{-1}_\gamma(x)) &\dots & m_{l_{n-1}} (-X_{n-1} g^{-1}\kappa^{-1}_\gamma(x)) \end{array} \right |.
\end{equation*}
Hence, $\det \Gamma(x,g)$ vanishes if, and only if, the fundamental vector fields $X_i(m)=X_i m$, $i=1,\dots, n-1$, do span the tangent spaces $T_m(p^{-1}(c)\cap V_\R)$ at each point $m=g^{-1} \kappa ^{-1}_\gamma(x)$, and if their span is not perpendicular to the hypersurface $\mklm{m\in V_\R: m_j=0}$. The latter condition is equivalent to $\grad p(m) \notin \mklm{m\in V_\R:m_j=0}$, which in turn is equivalent to the condition $\gd p /\gd m_j(m)\not=0$ for $m=g^{-1} \kappa^{-1}_\gamma(x)$. By construction, these conditions are fulfilled for $x \in \tilde U^\gamma$ and $g \in \supp c_\gamma(g)\subset G^\gamma_1$, as  a consequence of the properties (iii) and (iv) of the atlas $(\kappa^\gamma, U^\gamma)$. Thus, for arbitrary $\gamma$,
\begin{equation*}
  \Gamma(x,g) \text{ is invertible for } x\in \tilde U^\gamma, \, g \in \supp c_\gamma.
\end{equation*}
We now consider the extension of $\Gamma(x,g)$, as an endomorphism in $\C^1[\R^n_\xi]$, to the symmetric algebra ${\rm{S}}(\C^1[\R^n_\xi])\simeq \C[\R^n_\xi]$. Regarding the polynomials $\xi_1,\dots,\xi_n$ as a basis in $\C^1[\R^n_\xi]$, we obtain
\begin{equation*}
  \left ( 
  \begin{array}{c}  \Gamma \xi_1 \\ \vdots \\ \Gamma \xi_n \end{array} \right ) = (-i \tilde \phi _{-\xi,x}(g)) \left ( 
  \begin{array}{c} dR(X_1) \tilde \phi^\gamma_{\xi,x}(g) \\ \vdots \\ dR(Y)\tilde \phi^\gamma_{\xi,x}(g) \end{array} \right ),
\end{equation*}
 where $\Gamma(x,g)$ denotes the image of the basis vector $\xi_i$ under the endomorphism $\Gamma(x,g)$. Since, for $x \in \tilde U^\gamma$, $g \in \supp c_\gamma$, $\Gamma(x,g)$ is invertible, its extension to $ {\rm{S}}^N(\C^1[\R^n_\xi])$ is also an automorphism, so that every polynomial $\xi_{i_1} \otimes \dots \otimes \xi_{i_N}\equiv \xi_{i_1} \dots \xi_{i_n}$ can be written as a linear combination
 \begin{equation}
\label{24}
   \xi^\alpha =\sum _\beta \Lambda^\alpha_\beta (x,g) \Gamma \xi_{\beta_1} \dots \Gamma \xi_{\beta_{|\alpha|}}=-i\tilde \phi^\gamma_{\xi,x}(g) \sum \Lambda^\alpha_\beta(x,g) [dR(X_{\beta_1}) \tilde \phi^\gamma_{\xi,x}(g)]\dots,
 \end{equation}
where the $\Lambda^\alpha_\beta(x,g)$ are $\Cinft$ functions on $\tilde U^\gamma \times \supp c_\gamma$ which, in addition, are rational in $g$.
Now, we need  the following lemma.
\begin{lemma}
\label{lem:3}
  Under the identification $\Gamma \xi_k= -i \tilde \phi^\gamma_{\xi,x}(g) dR(X_k) \tilde \phi^\gamma_{\xi,x}(g) $ one has
  \begin{align}
\label{25}
\begin{split}
    i^r \tilde \phi^\gamma_{\xi,x}(g) \Gamma\xi_{j_1} \dots \Gamma\xi_{j_r}&= dR(X_{j_1} \dots X_{j_r})\tilde \phi^\gamma_{\xi,x}(g)\\&+ \sum_{s=1}^{r-1} \sum _{i_1,\dots, i_s} d ^{j_1,\dots, j_r}_{i_1,\dots, i_s} (x ,g) d R(X_{i_1} \dots X_{i_s}) \tilde \phi^\gamma_{\xi,x}(g),
\end{split}  
\end{align}
where the coefficients $ d ^{j_1,\dots, j_r}_{i_1,\dots, i_s} (x ,g) \in \Cinft(\tilde U^\gamma \times \supp c_\gamma)$ are at most of exponential growth in $g$, and independent of $\xi$.
\end{lemma}
\begin{proof}
  For $r=1$ one has $i \tilde \phi^\gamma_{\xi,x}(g) \Gamma \xi_k =d R(X_k) \tilde \phi^\gamma_{\xi,x}(g)$, and differentiating  the latter equation with respect to $X_j$ we obtain, by taking \eqref{24} into account,
  \begin{equation*}
    -\tilde \phi^\gamma_{\xi,x}(g) \Gamma \xi_j \Gamma \xi_k = dR(X_jX_k) \tilde \phi^\gamma_{\xi,x}(g)-\sum_{l,r} (dR (X_j) \Gamma_{kl}) (x,g) \Lambda ^l_r(x,g) dR (X_r) \tilde \phi^\gamma_{\xi,x}(g).
  \end{equation*}
Hence, the assertion of the lemma is correct for $r=1,2$. Now, assume that it holds for  $r\leq N$. Setting $r=N$ in equation \eqref{25}, and differentiating with respect to $X_k$, yields for the left hand side
\begin{align*}
  i^{N+1} \tilde \phi^\gamma_{\xi,x}(g) \Gamma \xi_k \Gamma\xi_{j_1} \dots \Gamma \xi_{j_N} &+ i^N \tilde \phi^\gamma_{\xi,x}(g) \Big ( \sum_{l,q} (dR(X_k) \Gamma_{j_1,l})(x,g) \Lambda_{l,q} (x,g) \Gamma \xi_q \Big ) \Gamma \xi _{j_2} \dots \Gamma \xi_{j_N}\\& + \dots.
\end{align*}
By assumption, we can apply \eqref{25} to the products $\Gamma \xi_q \Gamma\xi_{j_2} \dots\Gamma \xi_{j_N}$, and to all other products of fewer terms. This proves the lemma.
\end{proof}
Returning to the proof of Proposition \ref{prop:3}, as an immediate consequence of equations \eqref{24} and \eqref{25}, one computes
\begin{equation}
  \label{26}
(1+\xi^2)^N= \tilde \phi^\gamma_{-\xi,x}(g) \sum_{r=0}^{2N} \sum_{|\alpha| =r} b^N_\alpha(x,g) d R(X^\alpha) \tilde \phi^\gamma_{\xi,x}(g),
\end{equation}
where the coefficients $b^N_\alpha(x,g) \in \Cinft(\tilde U^\gamma \times \supp c_\gamma)$  are rational expressions in the matrix coefficients of $g$. Now, $\gd^\alpha_\xi \gd^\beta_x a ^\gamma_f(x,\xi)$ is a finite sum of terms of the form
\begin{equation*}
  \xi^\delta e^{-i(1,x_1,\dots,x_n) \cdot \xi} \int f(g) d_{\delta \beta}(x,g) \tilde \phi^\gamma_{\xi,x}(g)c_\gamma(g) d_G(g),
\end{equation*}
the functions $d_{\delta\beta}(x,g) \in \Cinft (\tilde U^\gamma \times \supp c_\gamma)$ being at most of exponential growth in $g$. Making use of equation \eqref{26}, and integrating by parts, we finally obtain for arbitrary $\alpha, \beta$ the estimate
\begin{equation*}
  |\gd ^\alpha_\xi \gd ^\beta _x \tilde a_f^\gamma (x,\xi) | \leq \frac 1 {(1+\xi^2)^N} C_{\alpha,\beta,\omega} \qquad x \in \omega,
\end{equation*}
where $\omega$ denotes an arbitrary  compact set  in $\tilde U^\gamma$ and $N=1,2,\dots$. This completes the proof of Proposition \ref{prop:3}.
\end{proof}
 In order to prove Theorem \ref{thm:3}, we need to examine the symbols $a_f^\gamma(x,\xi)$ and the kernels of the operators $A^\gamma_f$ more closely.
\begin{lemma}
  The kernel of $A_f^\gamma$ is given by the oscillatory integral
  \begin{equation}
\label{28}
    K_{A_f^\gamma}(x,y)=\int e^{i(x-y)\xi} a_f^\gamma(x,\xi) \, \dbar \xi,
  \end{equation}
where $a_f^\gamma(x,\xi)\in \Syms((\tilde U^\gamma_y \times \R^1_{\xi_1})\times \R^{n-1}_{\xi'})\cap \Sym^0((\tilde U^\gamma_y \times \R^{n-1}_{\xi'})\times \R^1_{\xi_1})$. It is a continuous function of $x\in \tilde U^\gamma$ with values in $\D'(\R^n_y)$.
\end{lemma}
\begin{proof}
  One computes $(\gd^\alpha_{\xi} a^\gamma_f)(x,\xi)=x_1^{\alpha_1} (\gd ^\alpha _\xi \tilde a^\gamma_f)(x,x_1\xi_1,\xi')$. By Proposition \ref{prop:3}, we therefore obtain the estimates
  \begin{equation*}
    |\gd ^\alpha_\xi\gd^\beta_y a^\gamma_f(x,\xi)| \leq C_{\alpha,\beta,K} |x_1|^{\alpha_1} \eklm{(x_1\xi_1,\xi')}^{l-|\alpha|}, \qquad l \in \R, \quad K \subset \subset \tilde U^\gamma.
  \end{equation*}
Hence, $a_f^\gamma(x,\xi)\in \Syms((\tilde U^\gamma_y \times \R^1_{\xi_1})\times \R^{n-1}_{\xi'})\cap \Sym^0((\tilde U^\gamma_y \times \R^{n-1}_{\xi'})\times \R^1_{\xi_1})$. In particular, $a^\gamma_f(x,\xi) \in \S(\R^{n-1}_{\xi'})$. Let $\chi\in \CT(\R)$ with $\chi(\xi_1)=1$ in a neighbourhood of $0$, and consider, for $u \in \CT(\R^n)$, the absolutely convergent integral
\begin{equation*}
  I_{x,\epsilon}(a_f^\gamma u)=\int \int e^{i(x-y)\xi} \chi(\epsilon \xi_1) a^\gamma_f(x,\xi) u(y) \, dy \dbar \xi.
\end{equation*}
 According to  \cite{shubin}, page 5, it can be regularized using integration by parts, so that, by the bounded convergence theorem, one may pass to the limit as $\epsilon \to 0$, in this way getting a distribution
\begin{equation*}
  u \mapsto I_{x}(a^\gamma_f u)= \lim_{\epsilon \to 0}I_{x,\epsilon}(a^\gamma_f u), \qquad u \in \CT(\R^n). 
\end{equation*}
As a consequence, \eqref{20} exists as oscillatory integral and is equal to $I_x(a^\gamma_f u)$, its kernel being given by \eqref{28}.
\end{proof}
We are now in position to prove Theorem \ref{thm:3}. 
\begin{proof}[Proof of Theorem \ref{thm:3}]
 According to  Proposition \ref{prop:3}, $\tilde a_f^\gamma(x,\xi) \in \Syms(\tilde U^\gamma \times \R^n_\xi)$, and it remains to show that $\tilde a_f^\gamma(x,\xi)$ satisfies the lacunary condition \eqref{V}.
Since, on $ {\tilde U^\gamma_+}$ resp. $ {\tilde U^\gamma_-}$,  $A^\gamma_f u$ only depends on the restriction of $u$ to ${\tilde U^\gamma_+}$ resp. $ {\tilde U^\gamma_-}$, we necessarily have
\begin{equation*}
  \supp K_{A^\gamma_f} \subset ({\tilde U^\gamma_+}\times {\tilde U^\gamma_+}) \cup ({\tilde U^\gamma_-}\times {\tilde U^\gamma_-}).
\end{equation*}
By the previous lemma, the kernel of $A^\gamma_f$ is given by the oscillatory integral
\eqref{28}. Consequently, 
\begin{equation*}
  \int e ^{i(x_1-y_1) \xi_1}  \tilde a_f^\gamma(x, x_1\xi_1, \xi') \, \dbar \xi_1,
\end{equation*}
which is a $\Cinft$ function for $x_1 \geq 0, y_1 <0$, resp.\ $x_1\leq 0, y>0$, must vanish then.  By making the substitution $t=y_1/x_1-1$ in case that $x_1\not=0$, we arrive at the condition
\begin{equation}
  M\tilde a_f^\gamma(x,\xi';t)=\int e^{-it\xi_1} \tilde a_f^\gamma(x,\xi) \, \dbar \xi_1 =0 \qquad \text{ for } t< -1, \, x \in \tilde U^\gamma.
\end{equation}
meaning that  $\tilde a_f^\gamma(x,\xi)$ satisifies the lacunary condition \eqref{V}. This completes the proof of Theorem \ref{thm:3}. 
\end{proof}

As a consequence of Theorem \ref{thm:3} we can write the kernel of $A_f^\gamma$, for $x_1\not= 0$, in the form
\begin{align}
\label{27}
\begin{split}
  K_{A_f^\gamma}(x,y) &= \int e ^{i(x-y)\cdot \xi} a^\gamma _f (x,\xi) \dbar \xi=\int e^{i(x-y) \cdot T_x^{-1} \xi} \tilde a_f^\gamma(x,\xi) |\det (T_x^{-1})'(\xi)| \dbar \xi\\
&=\frac 1 {|x_1|} \tilde A_f^\gamma(x,T_x^{-1}(x-y)), \qquad x_1 \not=0,
\end{split}
\end{align}
where $\tilde A_f^\gamma(x,y)$ denotes the inverse Fourier transform of $\tilde a_f^\gamma(x,\xi)$,
\begin{equation*}
  \tilde A_f^\gamma(x,y)= \int e ^{i y\cdot \xi} \tilde a_f ^\gamma(x,\xi) \, \dbar \xi.
\end{equation*}
Since, for $x \in \tilde U^\gamma$, $\tilde a_f^\gamma(x,\xi)$ is rapidly falling in $\xi$, 
 it follows that $\tilde A_f^\gamma(x,y) \in \S(\R^n_y)$,  the Fourier transform being an isomorphism on the Schwartz space.
Consider the restriction of the atlas $(\kappa_\gamma,U^\gamma)$ to $V_\R-S_\R$. By \eqref{27}, the restriction of the kernel of the operator $A_f^{\gamma}$ to the diagonal is given by 
\begin{equation*}
  k^\gamma_f(m)=K_{A^\gamma_f}(\kappa_\gamma(m),\kappa_\gamma(m))=\frac 1 {|p(m)|} \tilde A_f^{\gamma} (\kappa_\gamma(m),0),
\end{equation*}
and we denote the corresponding density on $V_\R-S_\R$ by $k_f$.
\begin{proposition}
\label{prop:4}
 Let $(G_C,\rho,V_\C)$ be a reductive prehomogeneous vector space whose singular set is an irreducible hypersurface, satisfying Assumptions $1$ and $2$. Let $V_\R-S_\R= V_1 \cup \dots \cup V_l$ be the decomposition of $V_\R-S_\R$ into its connected components, and set
  \begin{equation*}
    k^\gamma_f(m,s) = |p(m)| ^{s}  \tilde A_f^\gamma(\kappa_\gamma(m),0).
  \end{equation*}
Denote the corresponding density on $V_\R-S_\R$ by $k_{f,s}$. Then, for $\phi \in \CT(V_\R-S^{\rm{sing}}_\R)$, the integrals
\begin{equation}
  \int_{V_\R-S_\R} \phi \, k_{f,s} = \sum _{i=1}^l \int _{V_i} \phi \, k_{f,s}
\end{equation}
converge for $\Re s >0$, and can be continued analytically to meromorphic functions in $s \in \C$, satisfying functional equations of the form \eqref{K}.
\end{proposition}
\begin{proof}
  Let $\chi_\gamma$ be a partition of unity subordinated to the covering $\mklm{U^\gamma}$.  Let $dm$ denote Lebesgue measure in $V_\R$. Identifying the functions $k^\gamma_f(m,s)$ with the densities $k^\gamma_f(m,s)\d m$,  we obtain
    \begin{align*}
      \sum _{i=1}^l \int_{V_i} \phi \, k_{f,s} &= \sum_{i=1}^l \sum_\gamma \int_{V_i\cap U^\gamma} \phi(m) \chi_\gamma (m) k^\gamma_f(m,s) \d m\\&= \sum_{i=1}^l \sum_\gamma \int_{V_i\cap U^\gamma} \phi(m) \chi_\gamma (m) A^\gamma_f(\kappa_\gamma(m),0) |p(m)|^{s} d m .
    \end{align*}
Since $\phi(m) \chi_\gamma(m) \tilde A^\gamma_f(\kappa_\gamma(m),0) \in \CT(V_\R-S^{\rm{sing}}_\R)$, the assertion follows with Theorem \ref{thm:A}. Note that there exists a finite covering of $\supp \phi$ by the open sets $U^\gamma$, so that, in particular, the sum over $\gamma$ appearing in the last equation is finite.
\end{proof}

Let us illustrate this in the situation of Example \ref{ex:3}. By \eqref{22}, one has, in this case, the functional equations
\begin{equation*}
    \left ( \begin{array}{c} \int_{V_+^\ast} |p^\ast(\xi)|^{s-\frac n2} \hat \phi \ast \widehat{\chi_\gamma \tilde A^\gamma_f}(\xi) \d \xi \\\int_{V_-^\ast} |p^\ast(\xi)|^{s-\frac n2} \hat \phi \ast \widehat{\chi_\gamma\tilde A_f^\gamma}(\xi) \d \xi \end{array}\right ) = (2 \pi)^nC(s)  \left ( \begin{array}{c} \int_{V_+} k_{f,-s}^\gamma \chi_\gamma(m) \phi (m) \d m \\\int_{V_-} k_{f,-s}^\gamma \chi_\gamma(m) \phi (m) \d m \end{array}\right ),
  \end{equation*}
where we set $\tilde A_f^\gamma(m)=\tilde A_f^\gamma(\kappa_\gamma(m),0)$, and  $\phi \in \CT(V_\R-S_\R^{\rm{sing}})$.
\begin{remark} We close this section by noting that a trace of $\pi(f)$ can be defined in the situation of Theorem \ref{thm:3} and Proposition \ref{prop:4} by a regularization procedure analogous to the one in Proposition \ref{prop:2}. This is closely related to the definition of the $b$-trace in the theory of totally characteristic pseudodifferential operators.
\end{remark}

\section{Strongly elliptic differential operators and kernels of holomorphic semigroups}

As an application of the results of the previous sections, we study the holomorphic semigroup generated by a strongly elliptic differential operator associated with the left regular representation $(\pi,\Cvan(V_\R))$ of a prehomogeneous vector space $(G_\C,\rho,V_\C)$. Let $\pi$ be a continuous representation of a Lie group $G$ on a Banach  space $\B$, $\g$ the Lie-algebra of $G$, and
$X_1,\dots,X_d$ a basis of $\g$.  A differential operator
\begin{displaymath}
\Omega=\sum_{|\alpha| \leq l} (-i)^{|\alpha|} c_\alpha d\pi(X^\alpha)
\end{displaymath}
associated with the representation $\pi$ is called \emph{strongly elliptic}, if, for all  $\xi \in \R^d$, the relation $\Re \sum_{|\alpha|=l} c_\alpha \xi^\alpha \geq \kappa |\xi|^l$ is satisfied, where $\kappa>0$ is some positive number. By Langlands \cite{langlands}, the closure of a strongly elliptic operator generates a strongly continuous holomorphic semigroup of bounded operators on $\B$ given by
\begin{equation*}
S_t=\frac 1 {2\pi i} \int _\Gamma e^{\lambda t} ( \lambda \1+\overline{\Omega})^{-1} d\lambda,
\end{equation*}
where $\Gamma$ is a appropriate path in $\C$ coming from infinity and going to infinity, such that   $\lambda \notin \sigma(\overline{\Omega})$ for  $\lambda \in \Gamma$, and the integral converges uniformly with respect to the operator norm. Here $|\arg t|< \alpha$ for an appropriate $\alpha \in (0,\pi/2]$. The proof of this essentially relies on the verification of a criterium of Hille \cite{hille-phillips} within the theory of strongly continuous semigroups. The subgroup $S_t$ can be characterized by a convolution semigroup of complex measures $\mu_t$ on  $G$ according to  
\begin{equation*}
S_t=\int_G \pi(g) d\mu_t(g),
\end{equation*}
$\pi$ being measurable with respect to the measures $\mu_t$. The measures $\mu_t$ are absolutely continuous with respect to Haar measure $d_G$ on  $G$, and denoting by $K_t(g)\in L^1(G,d_G)$ the corresponding Radon-Nikodym derivative, one has
\begin{equation*}
S_t=\pi(K_t)=\int_G K_t(g) \pi(g) d_G(g).
\end{equation*}
The function $K_t(g)\in \L^1(G,d_G)$ is analytic in $t$ and $g$,
and given universally for all Banach representations. In the following, it will be called the \emph{Langlands kernel} of the holomorphic semigroup $S_t$. Though, in general, the Langlands kernel is not explicitely known, it satisfies the following $\L^1$- and $\L^\infty$-bounds. A detailed exposition  of these facts can be found in \cite{robinson}.
\begin{theorem}
Let $(L,\Cinft(G))$ be the left regular representation of $G$. Then, for each $\kappa \geq 0$, there exist constants $a,b,c>0$, and  $\omega\geq 0$ such that   
\begin{equation*}
\int_G |dL(X^\alpha)\gd^\beta_t K_t(g)| e^{\kappa |g|}\d_G(g) \leq a b^{|\alpha|}
c^\beta {|\alpha|}!\, \beta!(1+t^{-\beta-{|\alpha|}/l })e^{\omega t},
\end{equation*}
for all   $t>0$, $\beta=0,1,2,\dots$ and multiindices $\alpha$. Furthermore,
\begin{equation*}
|dL(X^\alpha)\gd^\beta_t K_t(g)|\leq a b^{|\alpha|}
c^\beta {|\alpha|}!\, \beta!(1+t^{-\beta-({|\alpha|}+d+1)/l })e^{\omega t}e^{-\kappa |g|},
\end{equation*}
for all $g \in G$, where $d=\dim \g$, and $l$ denotes the order of $\Omega$. In particular, $K_t \in \S(G)$.\end{theorem}

Let $(G_\C,\rho,V_\C)$ be a reductive prehomogeneous vector space defined over $\R$, and $(\pi,\Cvan(V_\R))$ the regular representation of $G_\R$ on  $V_\R$.  Assume that $\Omega$ is a strongly elliptic operator associated with $\pi$, and consider the corresponding holomorphic semigroup $S_t=\pi(K_t)$ with Langlands kernel $K_t$. It induces a continuous linear mapping 
\begin{equation}
\label{33}
S_t:\CT(\R^n) \longrightarrow \Cvan(\R^n) \subset \D'(\R^n),
\end{equation}
and since $K_t\in \S(G)$, the result of the previous sections concerning the Schwartz kernel of the operator \eqref{33} are applicable.  We will illustrate this for the classical heat kernel and the simplest prehomogeneous vector space, see Example \ref{ex:1}.

Thus, let $G_\C=\C^\ast$, $V_\C=\C$. In this case, $G_\R=\GL(1,\R)=\R^+_\ast$ acts multiplicatively on $V_\R=\R$. Let $a=1$ be a basis of  $\g=\R$. The action of  $d\pi(a)$ on the  G{\aa}rding subspace  $\Cvan(\R)_\infty$ corresponds on $\R$ to the vector field
\begin{displaymath}
d\pi(a)\phi(x)=\grad \phi(x)\cdot  (\frac d {dh} \e{-h} x\Big)_{\big\vert_{h=0}}=-x \frac d {dx} \phi(x).
\end{displaymath}
The Casimir operator $\Omega=\big( d\pi(a)\big )^2 \in \U$ of the considered   $G$-representation is therefore given by
\begin{equation*}
\Omega=x^2 \frac {d^2}{dx^2} + x \frac d {dx},
\end{equation*}
and constitutes a differential operator of Euler type. Let  $d_{G_\R}(x)=dx/x$ be Haar measure on $\GL(1,\R)$. Denoting by $S_t$ the holomorphic semigroup generated by $\overline{\Omega}$, and by $K_t(x)$ the corresponding  Langlands kernel, we get
\begin{equation*}
S_t(\phi)=\int_{\R^+_\ast} K_t(x) \pi(x)\phi  \frac {dx} x, \qquad \phi \in \Cvan(\R).
\end{equation*}
In order to compute $K_t(x)$ explicitely, we consider the Banach representation  $(\pi,\Cvan(\R^+_\ast))$, introducing on $\R^+_\ast$ the new coordinate   $x=e^r, \, r \in \R$. The Casimir operator is then given by $d^2/dr^2$, and the Langlands kernel coincides with the classical heat kernel  $K_t(r)=(2\pi t)^{-1/2} \exp({-r^2/4t})$. By transforming back, we get for $K_t(x)$ the universal expression 
\begin{equation*}
K_t(x)=\frac 1 {\sqrt{2\pi t}} \e{-(\log x)^2/4t}.
\end{equation*}
Now, 
\begin{displaymath}
(S_t\phi)(x)=\int_{\R^+_\ast} K_t(y) \phi(xy^{-1})  \frac {dy} y=\left
\{\begin{array}{c} \int\limits _0^\infty K_t(xy^{-1}) \phi(y) \frac {dy} y,
\qquad x >0,\\
\norm{K_t}_{L^1} \phi(0),\qquad \qquad \quad x=0,\\
 \int\limits ^0_{-\infty} K_t(xy^{-1}) \phi(y) \frac {dy}{|y|}, \quad\,\,\, x <0.\end{array}\right.
\end{displaymath}
If we therefore define 
\begin{equation*}
 S_t(x,y) =\left \{ \begin{array}{cc} K_t(xy^{-1})|y|^{-1} & \quad \text{for } xy^{-1}>0  \\ 0 & \quad \text{for } xy^{-1}\leq 0, \end{array} \right.
\end{equation*}
we get for $S_t$ the representation
\begin{equation*}
(S_t\phi)(x)= \left \{\begin{array}{cc} \int\limits_{-\infty}^\infty
S_t(x,y)\phi(y) dy, &
\quad x\not=0,\\
\norm{K_t}_{L^1} \phi(0),& \quad x=0.\end{array}\right.
\end{equation*}
One has $S_\R=\mklm{0}$, and Theorem \ref{thm:2} reads as follows.
\begin{proposition}
\label{prop:6}
The bounded linear operator $S_t:\CT(\R)\rightarrow\D'(\R)$ is given by the family of oscillatory integrals
\begin{equation*}
(S_t\phi)(x)=\int e^{i(x-y)\xi} s_t(x,\xi) \phi(y) dy \dbar \xi,
\end{equation*}
where 
\begin{equation*}
s_t(x,\xi)= e^{-ix\xi} \int_{\R^+_\ast} K_t(z) e^{\frac{i x\xi} z} \frac {dz} z \in \Sym^0_{1,0} ( \R_y \times \R_\xi), \qquad x \in \R;
\end{equation*}
for  $x\not=0$, one even has $s_t(x,\xi)\in \Syms (\R_y\times
\R_\xi)$. Its Schwartz kernel $\S_t\in \D'(\R\times \R)$ is given by the continuous family of oscillatory integrals 
\begin{equation*}
\R \ni x \mapsto \S_{t,x} = \int e^{i(x-\cdot) \xi} s_t(x,\xi) \dbar \xi.
\end{equation*}
Furthermore, $s_t(x,\xi) \in \Syms((\R_x-\mklm{0})\times \R_y\times \R_\xi)$, so that $S_t$, as a map  from   $\CT(\R)$ to   $\D'(\R-\mklm{0})$,    represents a Fourier integral operator. On $\R-\mklm{0}$,  $S_t\in\L^{-\infty} (\R-\mklm{0})$.
\end{proposition}
Since the Langlands kernel is explicitely known, it is illustrative to  give a direct proof in this case. Since $G$ is unimodular, one has
$$
s_t(x,\xi)= e^{-ix\xi} \frac 1 {2\pi t}\int\limits _0^\infty e^{-(\log y)^2/4t} e^{ixy\xi} \frac {d y} y = e^{-ix\xi} \hat f_t(-x\xi),
$$
where we defined the  function $f_t(y) =K_t(y) y^{-1}$, for $y>0$, $f_t(y) =0$ for $y\leq 0$. Clearly, $f_t \in \S(\R)$. Indeed, $K_t(y)y^m\in \Cinft(\R^+_\ast)$ for all  $m \in \Z$, and by the substitution $y=e^r$, we obtain
\begin{displaymath}
K_t(y)y^m=\frac 1 {\sqrt{2\pi t}} e^{-r^2/4t+mr}=\frac 1 {\sqrt{2\pi t}} e^{-(r/(2\sqrt t)-\sqrt t m)^2+m^2t};
\end{displaymath}
for $y\to 0$ and  $y\to\infty$, $K_t(y)y^m$ goes to zero together with all its derivatives, and it follows $|y^\alpha \gd^\beta_y f_t(y)|<\infty$ for arbitrary multiindices  $\alpha,\beta$. Hence we deduce $s_t(x,\xi)\in \Sym^1_{1,0}(\R_y\times \R_\xi)$ for arbitrary $x \in \R$, and  $s_t(x,\xi) \in \Syms(\R_y\times \R_\xi)$ for  $x\not=0$. Since    
\begin{align*}
\gd_\xi^\alpha \gd_x^\beta \hat f_t(-x\xi)=\sum\limits_{\alpha_1+\alpha_2=\alpha} \frac {\alpha !}{\alpha_1\alpha_2} \gd^{\beta+\alpha_1}_\chi \hat f_t(\chi)_{|\chi=-x\xi}(-x)^{\alpha_1} \beta\cdots (\beta-\alpha_2+1)(-\xi)^{\beta-\alpha_2},
\end{align*}
we finally obtain $s_t(x,\xi)\in \Syms(\R_x-\mklm{0}\times \R_y\times \R_\xi)$. The continuity of 
$x\mapsto S_{t,x}$ is a consequence of the fact that $\S_{t,x}(\phi)=(S_t\phi)(x)\in \Cvan(\R)$ is continuous. This proves Propostion \ref{prop:6}.

Now, by Lemma \ref{lem:2}, on $\R\times \R-\mklm{(0)}$, the Schwartz kernel $\S_t$ of  $S_t:\CT(\R)\rightarrow \D'(\R)$ is given by the function $S_t(x,y)$ introduced above. Indeed, for $x\not=0$,
\begin{align*}
\int e^{i(x-y)\xi} s_t(x,\xi) \dbar \xi=\int \int e^{i(z-y)\xi}
f_t(zx^{-1}) \frac 1 {|x|} \d z \dbar
\xi=f_t(yx^{-1})|x|^{-1}=S_t(x,y), 
\end{align*}
while, as oscillatory integral, 
\begin{displaymath}
\int e^{-iy\xi} s_t(0,\xi) \dbar \xi= \norm {K_t}_{\L^1} \, \delta_0,
\end{displaymath}
where $\delta_0$ denotes the  $\delta$-distribution mit support at zero. Furthermore, $S_t(x,y)\in\Lloc(\R^2)$. Following the discussion in Section \ref{sec:V}, we define     the auxiliary symbol
\begin{displaymath}
\tilde a_t(x,\xi)=e^{-i\xi} \hat f_t(-\xi) \in \Syms(\R_x \times \R_\xi),
\end{displaymath}
so that $s_t(x,\xi)=\tilde a_t(x,x\xi)$. The inverse Fourier transform       
\begin{equation*}
\tilde A_t(x,y)=\int e^{iy\xi} \tilde a_t(x,\xi) \dbar \xi = f_t(1-y)
\end{equation*}
is a Schwartz function, which vanishes for $y\geq 1$ together with all its derivatives. implying that  
$a_t(x,\xi)$ satisfies the lacunary condition \eqref{V}. For $x\not=0$ one computes
\begin{equation*}
S_t(x,y)=\int e^{i(x-y)\xi} \tilde a_t(x,x\xi) \dbar \xi=\frac 1 {|x|} \tilde A_t\big(x, \frac
{x-y} x\big ), \qquad x \not=0.
\end{equation*}
Furthermore, by the general theory of conormal distributions, $\tilde A_t(x,y) \in I^{-\infty}(\R\times
\R,\mklm{y=0})$, see   \cite{hoermanderIII}, Theorem $18.2.8$. In concordance with Proposition \ref{prop:2}, the restriction $s_t(x)=S_t(x,x)$ of $S_t$ to the diagonal defines a homogeneous distribution of degree $-1$  on $\R-\mklm{0}$, which has an extension to $\R$ given by \eqref{13}. Note that $\tilde A_t(x,0)=K_t(1)=(2\pi t)^{-1/2}$. Setting $s_{t,\zeta}(x)=(2\pi t)^{-1/2}|x|^\zeta$, \eqref{12} then yields in our situation the functional equations
\begin{equation*}
  \int ^\infty_{-\infty} s_{t,\zeta-1} \hat \phi(\xi) \d \xi=2 (2\pi)^{-\zeta} \Gamma(\zeta)  \cos \frac {\pi \zeta} 2 \int ^\infty _{-\infty} s_{t,-\zeta} (x) \phi (x) \d x, \qquad \phi \in \S(\R).
\end{equation*}

\providecommand{\bysame}{\leavevmode\hbox to3em{\hrulefill}\thinspace}
\providecommand{\MR}{\relax\ifhmode\unskip\space\fi MR }
\providecommand{\MRhref}[2]{%
  \href{http://www.ams.org/mathscinet-getitem?mr=#1}{#2}
}
\providecommand{\href}[2]{#2}

\end{document}